\documentclass[10pt]{article}

\usepackage{amssymb,amsmath}
\usepackage{multicol}
\usepackage{amsfonts}
  \usepackage{paralist}
\usepackage{enumerate}
\usepackage[latin1]{inputenc}

\usepackage[a4paper]{geometry}
\usepackage{graphicx}

\newtheorem{conj}{Conjecture}
\newtheorem{coro}{Corollary}

\newtheorem{defi}{Definition}
\newtheorem{lem}{Lemma}

\newtheorem{prop}{Proposition}
\newtheorem{theo}{Theorem}

\newtheorem{rem}{Remark}

\newtheorem{claim}{Claim}

\newenvironment{demo}
{\medbreak\noindent{\sc Proof :}} {\hfill $\square$ \medbreak}

\newenvironment{demof}[1]
{\medbreak\noindent{\sc Proof of {#1} :}}
{\hfill$\square$\medbreak}

\begin{document}

\newcommand\N{{\mathbb N}}
\newcommand\R{{\mathbb R}}

\title{\textbf{$\mathcal{C}^2$ surface diffeomorphisms have symbolic extensions}}
\author{David Burguet, \\
{\it CMLA-ENS Cachan}\\
{\it 61 avenue du président Wilson}\\ {\it 94235 Cachan Cedex France}
}

\date{}
\maketitle

\pagestyle{myheadings} \markboth{\normalsize\sc David
Burguet}{\normalsize\sc $\mathcal{C}^2$ surface diffeomorphisms have symbolic extensions}

\noindent \textbf{Abstract :} We prove that $\mathcal{C}^2$ surface diffeomorphisms have symbolic extensions, i.e. topological extensions which are subshifts over a finite alphabet. Following the strategy of T.Downarowicz and A.Maass \cite{Dow} we bound the local entropy of ergodic measures in terms of Lyapunov exponents. This is done by reparametrizing Bowen balls by contracting maps in a approach combining hyperbolic theory and Yomdin's theory.

\section{Introduction}

Given a dynamical system $(X,T)$ (i.e. a continuous map $T$ on a compact metrizable space $X$) one can try to encode it with a finite alphabet. More precisely we wonder if it admits a topological extension which is a subshift over a finite alphabet. Such an extension is called a symbolic extension. When a dynamical system has symbolic extensions  we are interested in minimizing their entropy. The topological symbolic extension entropy $h_{sex}(T)=\inf \{h_{top}(Y,S)$ :  $(Y,S)$  is a symbolic extension of $(X,T)\}$ estimates how the dynamical system $(X,T)$ differs from a symbolic extension from the point of view of entropy. These questions  lead to a deep theory which was developped mainly by M.Boyle and T.Downarowicz, who related the existence of symbolic extensions and their entropy with the convergence of the entropy of $(X,T)$ computed at finer and finer scales \cite{BD}. 

Dynamical systems with symbolic extensions have necessarily  finite topological entropy, but the converse is false : dynamical systems with finite  topological entropy may not have symbolic extensions. Nonetheless it was proved by M.Boyle, D.Fiebif, U.Fiebig \cite{bff} that asymptotically $h$-expansive dynamical systems with finite topological entropy  admit principal symbolic extensions, i.e. which preserve the entropy of invariant measures. Following Y.Yomdin \cite{Yoma},  J.Buzzi showed that $\mathcal{C}^{\infty}$ maps on a compact manifold are asymptotically $h$-expansive \cite{Buz}. In particular such maps admit principal symbolic extensions. On the other hand $\mathcal{C}^1$ maps without symbolic extensions  have been built in several works \cite{ND}, \cite{Asa}, \cite{bur}. T.Downarowicz and  A.Maass  have  recently proved  that $\mathcal{C}^r$ maps of the interval $f:[0,1]\rightarrow [0,1]$ with $1<r<+\infty$ have symbolic extensions \cite{Dow}. More precisely they showed that  $h_{sex}(f)\leq\frac{r\log\|f'\|_{\infty}}{r-1}$. The author built explicit examples  \cite{bur}  proving this upper bound is sharp. Similar $\mathcal{C}^r$ examples with large symbolic extension entropy have been previously built by T.Downarowicz and S.Newhouse for diffeomorphisms in higher dimension \cite{ND}. The results of T.Downarowicz and  A.Maass has been extended by the author in any dimension to nonuniformly entropy expanding maps (i.e. $\mathcal{C}^1$ maps whose ergodic measures with positive entropy have nonnegative Lyapunov exponent) of class $\mathcal{C}^r$ with $1<r<+\infty$ \cite{superbur}. The existence of symbolic extensions for general $\mathcal{C}^r$ maps with $1<r<+\infty$ is still an open question.

  We pointed out in the previous paragraph the known links between the regularity of the map and the existence of symbolic extensions. But symbolic extensions can also be obtained  by using geometrical arguments. For example it is well known that for uniformly hyperbolic dynamical systems, a Markov partition gives rise  to a principal symbolic extension \cite{Bow}. The elements of this partition are rectangles which are union of stable and unstable local manifolds. One can also just refer to the expansiveness of uniformly hyperbolic dynamical system and to the previously mentioned result of M.Boyle, D.Fiebif, U.Fiebig \cite{bff} to prove the existence of principal extensions. In fact $h$-expansiveness is also checked under a weaker assumption of hyperbolicity : partially hyperbolic dynamical system with a center bundle splitting into one dimensional subbundles \cite{cow} \cite{pac}.

 In the following we adapt such geometrical arguments to control the local dynamical complexity of  hyperbolic measures. More precisely we introduce 
  finite time rectangles which are pieces of finite time stable and  unstable manifolds. Under a condition of small oscillation of the derivative 
 these rectangles allows us to control the growth of the cardinality of separated sets. This condition can be ensured by the 
 $\mathcal{C}^2$ regularity of the map  as in Yomdin's theory. In this way this paper can be considered as a first attempt to combine the theory of hyperbolic dynamical systems and Yomdin's theory.

T.Downarowicz and S.Newhouse have conjectured in \cite{ND} that $\mathcal{C}^2$ diffeomorphisms on a compact manifold  have symbolic extensions. The following theorem answers affirmatively to this conjecture in the two dimensional case.  

\begin{theo}\label{intro}
Let $T:M\rightarrow M$ be a $\mathcal{C}^2$ surface diffeomorphism, then $T$ admits symbolic extensions and moreover 
$$h_{sex}(T)\leq h_{top}(T)+2R(T)$$
where $R(T)$ is the dynamical Lipschitz \footnote{$R(T)$ does not depend on the Riemannian metric $\|\|$ on $M$} constant of $T$, that is $R(T):=\lim_{n\rightarrow +\infty}\frac{\log^+\|DT^n\|}{n}$. 
\end{theo}

We recall in the second section the background of the theory of symbolic extensions. Then we state our main results and we reduce them to a theorem  of reparametrization of Bowen balls by contracting maps from the square in a similar (but finer) approach of Yomdin's $\mathcal{C}^2$ theory.  In the fourth section we introduce the finite time stable field which is a natural 
generalization  at finite time of the usual stable field. Assuming the oscillation of the derivative is small compared to the size of the derivative we prove the finite time stable field has bounded derivative. In the fifth section we define the new notion of finite time rectangle. Under the previous assumption on the oscillation of the derivative and some assumption of hyperbolicity we show these rectangles do not carry entropy. We also compare rectangles at successive times. The last section is devoted to the proof of the reparametrization of Bowen balls by such rectangles.

\section{Preliminaries}

In the following we denote $\mathcal{M}(X,T)$ the set of invariant Borel probability measures of the dynamical system $(X,T)$ and $\mathcal{M}_e(X,T)$ the subset of ergodic measures. We endow $\mathcal{M}(X,T)$ with the weak star topology. Since $X$ is a compact metric space, this topology is metrizable. We denote  $dist$ a metric on $\mathcal{M}(X,T)$. It is well known that $\mathcal{M}(X,T)$ is  compact and  convex and its extreme points  are exactly the ergodic measures. Moreover if $\mu\in \mathcal{M}(X,T)$ there exists an unique Borel probability measure $M_{\mu}$ on $\mathcal{M} (X,T)$ supported by $\mathcal{M}_e(X,T)$ such that  for all Borel subsets $B$ of $X$ we have $\mu(B)=\int\nu(B)dM_{\mu}(\nu)$. This is the so called ergodic decomposition of $\mu$. A Borel map $f:\mathcal{M}(X,T)\rightarrow \R$ is said to be harmonic if $f(\mu)=\int_{\mathcal{M}_e(X,T)}f(\nu)dM_{\mu}(\nu)$ for all $\mu\in \mathcal{M}(X,T)$. It is a well known fact that affine upper semicontinuous maps are harmonic.  In general\footnote{It may not be upper semicontinuous for $\mathcal{C}^r$ map for any $r\in \R^+$ \cite{Misss}. However it is for $\mathcal{C}^{\infty}$ maps \cite{Neww}.}  the measure theoretical entropy $h:\mathcal{M}(X,T)\rightarrow \mathbb{R}^+$ is not upper semicontinuous  but it is always harmonic \cite{W}.

If $f$ is a real Borel  map defined on $\mathcal{M}_e(X,T)$, the harmonic extension $\overline{f}$ of $f$ is the function defined on $\mathcal{M}(X,T)$ by :

$$\overline{f}(\mu):=\int_{\mathcal{M}_e(X,T)}f(\nu)dM_{\mu}(\nu)$$  
It is easily seen that $\overline{f}$ coincides with $f$ on $\mathcal{M}_e(X,T)$ and that $\overline{f}$ is harmonic.

\subsection{Entropy structure}
The measure theoretical entropy function can be computed in many ways as limits of a  nondecreasing sequence of  nonnegative functions defined on $\mathcal{M}(X,T)$ (with deacreasing sequence of partitions, formula of Brin-Katok,...). The entropy structures are such particular 
sequences whose convergence reflect the topological dynamic : they allow for example to compute the tail entropy \cite{Bur2} \cite{Dow2}, but also especially the symbolic 
extension entropy \cite{BD} \cite{Dow2}. 

We skip the formal definition of entropy structures. In the present paper we just use that Newhouse local entropy (See subsection \ref{newent}) is an entropy structure.

\subsection{Symbolic extension entropy function}
 A symbolic extension of $(X,T)$ is a subshift $(Y,S)$ of a full shift on a finite number of symbols, along with a continuous surjection $\pi: Y \rightarrow X$ such that $T\circ  \pi = \pi\circ S$. Given a symbolic extension $\pi:(Y,S)\rightarrow(X,T)$ we consider the extension entropy $h^{\pi}_{ext}:\mathcal{M}(X,T)\rightarrow \R^+$ defined for all $\mu\in \mathcal{M}(X,T)$ by : 
$$h^{\pi}_{ext}(\mu)=\sup_{\pi^*\nu=\mu}h(\nu)$$

Then the symbolic entropy function, $h_{sex}:\mathcal{M}(X,T)\rightarrow \R^+$,  is :
$$h_{sex}=\inf h_{ext}^{\pi}$$

where the infimum holds over all the symbolic extensions of $(X,T)$. By convention, if $(X,T)$ does not admit any symbolic extension we simply put 
$h_{sex}\equiv+\infty$. Recall we define in the introduction the topological symbolic extension entropy $h_{sex}(T)$ as the infimum of the topological entropy of the symbolic extensions of $(X,T)$ (as previously we put $h_{sex}(T)=+\infty$ when there are no such extensions). M.Boyle and T.Downarowicz proved that these two quantities are related by the following variational principle :

\begin{eqnarray}\label{pvsex}
h_{sex}(T)&=&\sup_{\mu\in \mathcal{M}(X,T)}h_{sex}(\mu)
\end{eqnarray}

\subsection{The Estimate Theorem}

One of the main tools introduced in \cite{Dow} is the so called Estimate Theorem. We can roughly resume the statement as follows : in order to estimate the symbolic entropy function one only needs to bound the local entropy of an ergodic measure  near an invariant one by the difference of the values of some upper semicontinuous affine function on $\mathcal{M}(X,T)$ at these two measures.

\begin{theo}[Downarowicz, Maass]\cite{Dow}\label{theopass}
Let $(X,T)$ be a dynamical system with finite topological entropy. Let $\mathcal{H}=(h_k)_k$ be an entropy structure. Let $g$ be an upper semicontinuous affine function on $\mathcal{M}(X,T)$ which dominates the entropy function at all ergodic measures and such that for every $\gamma>0$ and $\mu\in \mathcal{M}(X,T)$ there exists  $k_{\mu}\in\N$ and $\delta_{\mu}>0$ such that for every ergodic measure $\nu$ satisfying $dist(\nu,\mu)<\delta_{\mu}$ it holds that :

$$(h-h_{k_{\mu}})(\nu)\leq g(\mu)-g(\nu)+\gamma$$
Then, for every invariant measure $\mu$ on $X$, 
   $$h_{sex}(\mu)\leq h(\mu)+g(\mu)$$

\end{theo}

It is worth noting that the Estimate Theorem holds in a pure topological setting. This result is in fact a consequence of the major Theorem of Symbolic Extension \cite{BD} which allows to compute the symbolic extension entropy function from the properties of convergence of any entropy structure.

\subsection{Newhouse local entropy}\label{newent}
Let us first recall some usual notions relating to the entropy of dynamical systems (we refer to \cite{W} for a general introduction to entropy). Consider a continuous map $T:X\rightarrow X$ with $(X,d)$  a compact metric space. Let $n\in \N$ and $\delta>0$. A subset $F$ of $X$ is called a $(n,\delta)$ separated set when for all $x,y\in F$ there exists $0\leq k< n$ such that $d(f^kx,f^ky)\geq\delta$. Let $Y$ be a subset of $X$. A subset $F$ of $Y$ is called a $(n,\delta)$ spanning  set of $Y$ when for all $y\in Y$ there exists $z\in F$ such that $d(f^k(x),f^k(y))<\delta$ for all $0\leq k< n$. Given a point $x\in X$ we denote by $B(x,n,\delta)$ the Bowen ball centered at $x$ of radius $\delta$ and length $n$ :
$$B(x,n,\delta):=\{y\in X, \ d(T^ky,T^kx)<\delta \ \text{for} \ k=0,...,n-1\}$$
 
We recall now the "Newhouse local entropy". Let $x\in X$, $\epsilon>0$, $\delta>0$, $n\in \N$ and $F\subset X$ a Borel set, we define :

$$H(n,\delta|x,F,\epsilon):=\log \max\left\{\sharp E \ : \  E\subset F\bigcap B(x,n,\epsilon) \ \text{and} \  E \ \text{is a}  \ (n,\delta) \ \text{separated set}\right\} $$

$$H(n,\delta|F,\epsilon):=\sup_{x\in F} H(n,\delta|x,F,\epsilon)$$

$$h(\delta|F,\epsilon):=\limsup_{n\rightarrow +\infty}\frac{1}{n}H(n,\delta|F,\epsilon)$$

$$h(X|F,\epsilon):=\lim_{\delta\rightarrow 0}h(\delta|F,\epsilon)$$

 Then for any ergodic measure $\nu$ we put :
$$h^{New}(X|\nu,\epsilon):=\lim_{\sigma\rightarrow 1}\inf_{\nu(F)>\sigma}h(X|F,\epsilon)$$

Given a nonincreasing sequence $(\epsilon_k)_{k\in \N}$ converging to $0$, we consider the sequence
$\mathcal{H}^{New}=(h_k^{New})_{k\in\N}$ where $h_k^{New}:=\overline{h-h^{New}(X|.,\epsilon_k)}$ is the harmonic extension of $h-h^{New}(X|.,\epsilon_k)$ for all integers $k$. T.Downarowicz proved this sequence defines an entropy structure \cite{Dow2} for homeomorphisms and the author extends the result in the noninvertible case \cite{Bur}. In particular $h^{New}(X|.,\epsilon_k)$ converges pointwise to zero when $k$ goes to infinity. \\

We will use the following technical inequality :

\begin{lem}\label{tecnew}
For all integers $k>0$ and for all ergodic measures $\nu$ :
$$h_{T}^{New}(X|\nu,\epsilon)\leq \frac{h_{T^{k}}^{New}(X|\nu,\epsilon)}{k}$$
\end{lem}

\begin{demo}
Clearly we have the inclusion $B_T(x,nk,\epsilon)\subset B_{T^k}(x,n,\epsilon)$. Moreover for all $\delta>0$ there exists $\delta'>0$ such that any $(nk,\delta)$ separated set for $T$ is $(n,\delta')$ separated for $T^k$. The conclusion of the lemma follows easily from these two facts.
\end{demo}

\subsection{Lyapunov exponents for surface diffeomorphisms}

Let $(M,\|\|)$ be a compact Riemannian surface and let $T:M\rightarrow M$ be a diffeomorphism and $\nu$ be an ergodic $T$-invariant measure.  We denote $\|D_xT\|=\sup_{v\in T_xM-\{0\}}\frac{\|D_xT(v)\|}{\|v\|}$ the induced norm of the differential $D_xT$ of $T$ at $x$  and  $\|DT\|=\sup_{x\in M}\|D_xT\|$ the supremum norm of the differential of $T$. According to Oseledet's theorem \cite{Osc}, there exist two real numbers  $\chi^+(\nu)\geq 
\chi^-(\nu)$, a measurable splitting of the tangent bundle $TM=E_1\bigoplus E_2$ into two invariant subbundles $E_1$ and $E_2$ and a Borel set $F$ 
with $\nu(F)=1$ such that for all $x\in F$ and all $(v_1,v_2)\in E_1\times E_2$ :
\begin{eqnarray*}
  \lim_{|n|\rightarrow +\infty}\frac{1}{n}\log\|D_xT^n(v_1)\|&=&\chi^+(\nu)\\
  \lim_{|n|\rightarrow +\infty}\frac{1}{n}\log\|D_xT^n(v_2)\|&=&\chi^{-}(\nu)
\end{eqnarray*}

Remark that $\chi^+(\nu)=\lim_{n\rightarrow +\infty}\frac{1}{n}\log\|D_xT^n\|$ and $\chi^-(\nu)=\lim_{n\rightarrow -\infty}\frac{1}{n}\log\|D_xT^{n}\|$ for all $x\in F$. The real numbers $\chi^+(\nu)$ and $\chi^-(\nu)$ are the well-known Lyapunov exponents of $\nu$ (we use also the notations $\chi^+(\nu,T)$ and $\chi^-(\nu,T)$ to avoid ambiguities). We denote $\chi^+_0(\nu)=\max(\chi^{+}(\nu),0)$ and $\chi^-_{0}(\nu)=\min(\chi^-(\nu),0)$. According to the subadditive ergodic theorem, we have  $\chi^+_0(\nu)=\inf_{n\in\N}\frac{1}{n}\int_M\log^+\|D_xT^n\|d\nu(x)=\lim_{n\rightarrow +\infty}\frac{1}{n}\int_M\log^+\|D_xT^n\|d\nu(x)$ and  similarly $-\chi^-_0(\nu)=\inf_{n\in\N}\frac{1}{n}\int_M\log^+\|D_xT^{-n}\|d\nu(x)=\lim_{n\rightarrow +\infty}\frac{1}{n}\int_M\log^+\|D_xT^{-n}\|d\nu(x)$. Another application of the subadditive ergodic theorem leads to the following lemma : 

\begin{lem}\label{lyap}
Let $\nu$ be an ergodic measure, then
$\frac{1}{n}\log^+\|(D_xT^{n})^{-1}\|$ converges almost everywhere to $-\chi^-_0(\nu)$ when $n$ goes to $+\infty$.
\end{lem}

\begin{demo}
By the subadditive ergodic theorem, the sequence $(\frac{1}{n}\log^+\|(D_xT^{n})^{-1}\|)_{n\in\N}$ converges $\nu$ almost everywhere to 
$\inf_{n\in\N}\frac{\int\log^+\|(D_xT^{n})^{-1}\|d\nu(x)}{n}$ when $n$ goes to $+\infty$. Then we have by invariance of $\nu$ :
\begin{eqnarray*}
 \int\log^+\|(D_xT^{n})^{-1}\|d\nu(x)&=&\int\log^+\|D_{T^nx}T^{-n}\|d\nu(x) \\
 &=&\int\log^+\|D_{x}T^{-n}\|d\nu(x) 
 \end{eqnarray*}
\end{demo}

We prove now elementarily that the harmonic extension of $\chi^+_0$ is upper semicontinuous. The following lemma is in fact valid in any dimension. In the one dimensional case it was proved by T.Downarowicz and A.Maass by using a clever argument of convexity (See Fact 2.5 of \cite{Dow}).   

\begin{lem}\label{thes}
For all $\mu\in \mathcal{M}(M,T)$, we have :

$$\overline{\chi^+_0}(\mu)=\inf_{n\in\N}\frac{1}{n}\int_M\log^+\|D_xT^n\|d\mu(x)$$

 In particular $\overline{\chi^+_0}:\mathcal{M}(M,T)\rightarrow \R^+$ is upper semicontinuous.  
\end{lem}

\begin{demo}
For all integers $n>0$ we consider the function $f_n:\mathcal{M}(M,T)\rightarrow \R^+$ defined by :
$$\forall \mu \in \mathcal{M}(M,T), \ f_n(\mu)=\int \log^+\|D_xT^n\|d\mu(x)$$
This function is clearly continuous and affine, and therefore harmonic. Also  $(f_n(\mu))_{n\in\N}$ is a subadditive sequence for all $\mu\in \mathcal{M}(M,T)$. 

We already observe that $\chi^+(\nu)=\lim_{n\rightarrow +\infty}\frac{f_n(\nu)}{n}$ for all ergodic measures $\nu$.
Consider now a general measure $\mu\in \mathcal{M}(M,T)$. We have :

\begin{eqnarray*}
\overline{\chi^+_0}(\mu)&:=&\int_{\mathcal{M}_e(M,T)}\chi^+_0(\nu)dM_{\mu}(\nu)\\
 &=&\int_{\mathcal{M}_e(M,T)}\lim_{n\rightarrow +\infty}\frac{f_n(\nu)}{n}dM_{\mu}(\nu)
\end{eqnarray*}

Obvously $f_n(\nu)\leq \log^+\|DT\|$ for all ergodic measures $\nu$. Therefore by applying the theorem of dominated convergence we get :

\begin{eqnarray*}
\overline{\chi^+_0}(\mu)&=&\lim_{n\rightarrow +\infty}\int_{\mathcal{M}_e(M,T)}\frac{f_n(\nu)}{n}dM_{\mu}(\nu)
\end{eqnarray*}

and by harmonicity of $f_n$ :

\begin{eqnarray*}
\overline{\chi^+_0}(\mu)&=&\lim_{n\rightarrow +\infty}\frac{f_n(\mu)}{n}
\end{eqnarray*}

But the sequence $(f_n(\mu))_{n\in \N}$ is subadditive so that :

\begin{eqnarray*}
\overline{\chi^+_0}(\mu)&=&\inf_{n\in\N}\frac{f_n(\mu)}{n}
\end{eqnarray*}

We conclude that $\overline{\chi^+_0}$ is an upper semicontinuous function as an infimum of a family of continuous functions. 
\end{demo}

In the following we are interesting in the entropy of ergodic measures. Recall the  Ruelle-Margulis inequality states that for a $\mathcal{C}^1$ map $T:M\rightarrow M$ on a compact manifold $M$ the entropy $h_T(\nu)$ of an ergodic measure $\nu$ is bounded from above by the sum of its positive Lyapunov exponents. When $T$ is a surface diffeomorphism it is easily seen  by applying the Ruelle-Margulis inequality to $T$
 and its inverse $T^{-1}$ and by using the equality $h_T(\nu)=h_{T^{-1}}(\nu)$ that  any ergodic measure $\nu\in \mathcal{M}(M,T)$ with non zero entropy has exactly one positive and one negative Lyapunov exponent, i.e., with the previous notations,  $\chi^{+}(\nu)>0$ and $\chi^-(\nu)<0$ and moreover 

$$h(\nu)\leq \min(\chi^+(\nu),-\chi^-(\nu))$$

\section{Statements}

\begin{theo}\label{super}(Main Theorem)
Let $T:M\rightarrow M$ be a $\mathcal{C}^2$ surface diffeomorphism. Then for all $\mu\in \mathcal{M}(M,T)$, 

$$h_{sex}(\mu)\leq h(\mu)+2\overline{\chi^+_0}(\mu)$$

\end{theo}

The Main Theorem follows easily from the following theorem by applying the Estimate Theorem to the upper semicontinuous affine function $g_0=2\overline{\chi^+_0}$ (which dominates the entropy function at all ergodic measures according to Ruelle-Margulis inequality) and to the entropy structure $\mathcal{H}^{New}$. Remark  that Theorem \ref{intro} stated in the introduction is the topological version of the Main Theorem : it is deduced by the usual variational principle for the entropy and the variational principle for the symbolic extension entropy (Equation (\ref{pvsex})).

\begin{theo}\label{gg}
Let $T:M\rightarrow M$ be a $\mathcal{C}^2$ surface diffeomorphism. Let $\mu$ be an invariant measure and fix some
$\gamma> 0$. Then there exist $\delta_{\mu}>0$ and $\epsilon_{\mu}> 0$ such that for
every ergodic measure $\nu$ with $dist(\nu,\mu) < \delta_{\mu}$ it holds that :

\begin{eqnarray}\label{ggeq}
h^{New}(M|\nu,\epsilon_{\mu})\leq 2\overline{\chi^+_0}(\mu)-2\chi^+_0(\nu)+\gamma \nonumber
\end{eqnarray}

\end{theo}

In the general case (any dimension, any intermediate regularity, noninvertible maps) we conjecture that Theorem \ref{gg} can be extended in the following way : 

\begin{conj}\label{mainconj}
Let $r>1$. Let $M$ be a compact manifold of dimension $d$ and  $T:M\rightarrow M$  a $\mathcal{C}^r$ map\footnote{for $r\notin \N$ we mean that $T\in \mathcal{C}^{[r]}$ and $D^{[r]}T$ is $r-[r]$-Hölder}. Let $\mu$ be an invariant measure and fix some $\gamma> 0$. Then there exist $\delta_{\mu}>$ and $\epsilon_{\mu}> 0$ such that for
every ergodic measure $\nu$ with $dist(\nu,\mu) < \delta_{\mu}$ it holds that :

\begin{eqnarray}
h^{New}(M|\nu,\epsilon_{\mu})\leq \frac{\overline{\Sigma^+\chi}(\mu)-\Sigma^+\chi(\nu)}{r-1}+\gamma \nonumber
\end{eqnarray}
where $\Sigma^+\chi$ denotes the sum of the positive Lyapunov exponents. 
\end{conj}

\begin{rem}
By adapting the proof of Lemma \ref{thes} one easily shows that $$\overline{\Sigma^+\chi}(\mu)=\inf_{n\in\N}\frac{1}{n}\int \max_{k=1,...,d}\log^+\|\Lambda^kD_xT^n\|_kd\mu(x)$$ and therefore $\overline{\Sigma^+\chi}$ is also an upper semicontinuous function on $\mathcal{M}(M,T)$.
\end{rem}

This conjecture was proved by T.Downarowicz and A.Maass in dimension one \cite{Dow}. In \cite{superbur} we prove the conjecture in any dimension up to a factor $d$ and an additional term corresponding to the sum of the negative Lyapunov exponents of $\nu$  (then to apply the Estimate Theorem we are reduced to the class of nonuniformly entropy expanding maps because of this remaining term). Theorem \ref{gg} corresponds to the case of  $\mathcal{C}^2$ surface diffeomorphisms up to a factor $2$. 

\begin{rem}
By the Passage theorem of T.Downarowicz and A.Maass \cite{Dow} the inequality in Theorem \ref{gg} holds in fact for general invariant measures $\nu$ close to $\mu$. Then 
following the proof of the Estimate Theorem we get that $u_1(\mu)\leq \frac{2\overline{\chi^+_0}(\mu)}{3}$ for all invariant measures $\mu$. Taking the supremum over all invariant measures $\mu$ we have according to the tail variational principle  the following bound on the tail\footnote{We refer to \cite{Dow2} and to \cite{Bur2} for the definitions of $u_1$ and $h_{tail}(T)$ and the tail variational principle.} entropy : $h_{tail}(T)\leq \frac{2R(T)}{3}$. By using Yomdin's theory, J.Buzzi \cite{Buz} proved that $h_{tail}(S)\leq \frac{d}{r}R(S)$ when $S$ is a  $\mathcal{C}^r$ map with $r\geq 1$ on a compact manifold of dimension $d$. Here we consider $\mathcal{C}^2$ surface diffeomorphism ($d=r=2$) so that our upper bound is better than Buzzi's one. It is known \cite{Buz} that Buzzi's inequality is sharp for noninvertible maps but it is reasonnable to think that only the expanding directions are involved in the creation of tail entropy as in the Ruelle-Margulis inequality for the entropy (see also the previous conjecture) and thus that $h_{tail}(T)\leq \frac{R(T)}{2}$ for a $\mathcal{C}^2$ surface diffeomorphism $T$.
\end{rem}

We reduce now Theorem \ref{gg} to a result of reparametrization of Bowen's balls by contracting maps as in Yomdin's theory. 
Let us first introduce the notion of finite time hyperbolic sets. We fix a Riemannian metric $\|\|$ on the surface $M$ and we denote by $d$ the induced distance on $M$.

\begin{defi}
For any $\chi^+>0>\chi^-$, $\min(\chi^+,-\chi^-)>\gamma>0$ and $C>1$, we denote for all integers $n$ :
 
\begin{multline*}
\mathcal{H}^n_T(\chi^+,\chi^-,\gamma,C):=\{x\in M \ : \ \forall 1\leq k\leq n, \ C^{-1}e^{(\chi^+-\gamma)k}\leq \|D_xT^k\|\leq Ce^{(\chi^++\gamma)k}, \\  
C^{-1}e^{(-\chi^--\gamma)k}\leq \|D_{T^kx}T^{-k}\|\leq Ce^{(-\chi^-+\gamma)k}\}
\end{multline*}

\end{defi}
This set captures the hyperbolicity with logarithmic expansion $\chi^+$ and contraction $\chi^-$ with error $\gamma$ at finite time.
The following  proposition bounds the local dynamical complexity of finite time hyperbolic sets. We denote $H:[1,+\infty[\rightarrow \R$ the function defined by $H(t)=-\frac{1}{t}\log(\frac{1}{t})-(1-\frac{1}{t})\log(1-\frac{1}{t})$. Moreover $[x]$ is the integer part of $x$ if  $x>0$ and zero if not. Finally we say that a map from the unit square $[0,1]^2$ to a surface $M$ is $\mathcal{C}^1$ if it can be extended in a $\mathcal{C}^1$ map on an open neighborhood of $[0,1]^2$. 
 
\begin{prop}\label{main}
Let $T:M\rightarrow M$ be a surface diffeomorphism and let $\chi^+>0>\chi^-$, $\frac{\min(\chi^+,-\chi^-,1)}{3}>\gamma>0$ and $C>1$. Then  there exist $\epsilon>0$  depending only on $\|D^2T\|$, $\|D^2T^{-1}\|$, $\|DT\|$ and $\|DT^{-1}\|$,  a real number $D$ depending only on  $\chi^+,\chi^-,\gamma,C,\|DT\|$ and $\|DT^{-1}\|$ and a universal constant $A>0$  with the 
following properties. For all  $x\in M$ and for all positive integers $n$ there exists a family $\mathcal{F}_n$ of $\mathcal{C}^1$ maps from $[0,1]^{2}$ to $M$ 
such that :

\begin{enumerate}[(i)]
\item $$\forall \phi\in \mathcal{F}_n \ \forall 0 \leq l\leq n, \ \|D(T^{l}\circ \phi)\|\leq 1$$   
\item $$\mathcal{H}^n_T(\chi^+,\chi^-,\gamma,C)\cap B(x,n+1,\epsilon)\subset \bigcup_{\phi\in \mathcal{F}_n}\phi([0,1]^2)$$ 
\item $$\log\sharp \mathcal{F}_n \leq  \left(2+H([\lambda^+_n(x,T)-\chi^+]+3)\right)\left(\lambda^+_n(x,T)-\chi^+\right)n+An+D$$

 with $\lambda^+_n(x,T):=\frac{1}{n}\sum_{l=0}^{n-1}\log^+\|D_{T^lx}T\|$.
\end{enumerate}
\end{prop}

When the size $\epsilon$ of the Bowen ball at $x$ is small, $\lambda^+_n(y,T)$ is close to $ \lambda^+_n(x,T)$ for all $y$ in the Bowen ball. Then the term $\lambda^+_n(x,T)-\chi^+$ is (up to a  small error term) the logarithm of the defect of multiplicativity of the norm of the composition $D_yT^n=D_{T^{n-1}y}T\circ...\circ D_yT$ when $y$ belongs also to the finite hyperbolic set $\mathcal{H}^n_T(\chi^+,\chi^-,\gamma,C)$.\\ 

%Assume furthermore there exists a point $z$ in the Bowen ball such that  $\lambda^+_n(z,T)$ is closed to $\frac{1}{n}\log\|D_zT^n\|$, then the quantity $\lambda^+_n(x,T)-\chi^+$ estimates how the small oscillations of the norm of the differential along the orbit of $x$ accumulate to create a gap between the  maximal logarithmic expansion of $DT^n_{/B(x,n,\epsilon)}$ which is  $\frac{1}{n}\log\|D_zT^n\|\simeq \lambda^+_n(x,T)$ and the logarithmic expansion at points of $\mathcal{H}^n_T(\chi^+,\chi^-,\gamma,C)\cap B(x,n,\epsilon)$ which is nearly $\chi^+$.\\

We  deduce now Theorem \ref{gg} from the above statement. In the proof the terms $\lambda^+_n(x,T)$ for typical $\nu$ points $x$ and $\chi^+$ will be respectively related with the Lyapunov exponents of $\mu$ and $\nu$ where $\nu$ is an ergodic measure near an invariant measure $\mu$. Moreover the quantity  $H([\lambda^+_n(x,T)-\chi^+]+3)$ will be negligible. 

\begin{demof}{Theorem \ref{gg} assuming Proposition \ref{main}}

Let  $\mu\in \mathcal{M}(M,T)$. By Lemma \ref{thes} we choose $k_{\mu}\in \N$ s.t. 
\begin{eqnarray*}
\overline{\chi^+_0}(\mu,T)=\inf_{n\in \N}\frac{\int\log^+\|D_xT^n\|d\mu(x)}{n}&\geq &\frac{\int\log^+\|D_xT^{k_{\mu}}\|d\mu(x)}{k_{\mu}}-\gamma \label{ffirst}
\end{eqnarray*}

One can also assume $k_{\mu}$ large enough s.t. $2H(k_{\mu}\gamma)R(T)<\gamma$ and $\frac{A}{k_{\mu}}<\gamma$.

Then by continuity of $\mu\mapsto \int\log^+\|D_xT^{k_{\mu}}\|d\mu(x)$ and by upper semicontinuity of $\overline{\chi_0^{+}}$  one can choose the parameter $\delta_{\mu}>0$ such that for all ergodic measures $\nu$ with $dist(\nu,\mu)<\delta_{\mu}$ we have : 
 
 $$\left| \int \log^+\|D_{x}T^{k_{\mu}}\|d\nu(x)-\int\log^+\|D_xT^{k_{\mu}}\|d\mu(x) \right| < \gamma$$
 
 $$\overline{\chi_0^+}(\mu,T)>\chi_0^+(\nu,T)-\gamma $$

We fix some ergodic measure $\nu$ with $dist(\mu,\nu)<\delta_{\mu}$. One can assume $h(\nu)>3\gamma$. By Ruelle-Margulis inequality we have therefore  $\chi^+(\nu)=\chi^+_0(\nu)>3\gamma$ and $\chi^-(\nu)=\chi^-_0(\nu)<-3\gamma$.

According to the subadditive ergodic theorem (see Lemma \ref{lyap}), there exists for any $0<\sigma<1$ a Borel set $F_\sigma$  of $\nu$ measure larger than $\sigma$ s.t. $\exists n_{\sigma} \ \forall n>n_{\sigma}$,  $\left| \chi^+(\nu)-\frac{1}{k_{\mu}n}\log^+\|D_{x}T^{k_{\mu}n}\| \right| <\gamma$ and $\left| \chi^-(\nu)-\frac{1}{k_{\mu}n}\log^+\|D_{T^{k_{\mu}n}x}T^{-k_{\mu}n}\| \right| <\gamma$. 
Remark that these inequalities  can be rewritten as :

$$\exists C>1 \ \forall n\in \N, \ F_{\sigma}\subset \mathcal{H}_{T^{k_{\mu}}}^n(\chi^+(\nu),\chi^-(\nu),\gamma,C)$$ 

One can also assume by the ergodic theorem  that $(\lambda_n^+(x,T^{k_{\mu}}))_{n\in \N}$ are converging uniformly in $x\in F_{\sigma}$ to $ \int \log^+\|D_yT^{k_{\mu}}\|d\nu(y)$. \\

 We apply now Proposition \ref{main} to $T^{k_{\mu}}$  and to a given point $x\in F_{\sigma}$ : there exists $\epsilon_{\mu}>0$  depending only on $\|D^2T^{k_{\mu}}\|$, $\|D^2T^{-k_{\mu}}\|$, $\|DT^{k_{\mu}}\|$ and $\|DT^{k_{\mu}}\|$, a real number $D$ depending only on $\chi^+(\nu),\chi^-(\nu),\gamma,C,\|DT^{k_{\mu}}\|$ and $\|DT^{-k_{\mu}}\|$ and families $(\mathcal{F}_n)_{n\in\N}$ of $\mathcal{C}^1$ maps from $[0,1]^{2}$ to $M$ satisfying the properties (i),(ii),(iii) of Proposition \ref{main}.

 For all $0<\delta<1$ let $E_{\delta}\subset[0,1]^2$ denote the subset of points of the square of the form $k\delta$ with $k\in\N^2$. Consider  $\phi\in \mathcal{F}_n$. It follows from the first item (i) that $\phi(E_{\delta})$ is a 
$(n+1,\delta)$ spanning set of $\phi([0,1]^2)$ for $T^{k_{\mu}}$. Indeed given a point $y=\phi(z)$ with $z\in [0,1]^2$ there exists $t\in E_{\delta}$ 
such that $\|z-t\|<\delta$ and then by using the first item (i), we get $d\left(T^{lk_{\mu}}y,T^{lk_{\mu}}(\phi(t))\right)<\delta$ for $l=0,...,n$, that is $y\in B_{T^{k_{\mu}}}(\phi(t),n+1,\delta)$. Since two points lying in the same Bowen ball of radius $\delta$ and length $n+1$ are not $(n+1,2\delta)$ separated, we get :    
$$\max \left\{\sharp E \ : \ E\subset F_{\sigma}\cap \phi([0,1]^2)  \  \text{and} \ E \ \text{is a} \ (n+1,2\delta) \  \text{separated  set for} \  T^{k_{\mu}} \right\} \leq \frac{1}{\delta^2}$$
and then 
$$\max \left\{\sharp E \ : \ E\subset F_{\sigma}\cap B(x,n+1,\epsilon)  \  \text{and} \ E \ \text{is a} \ (n+1,2\delta) \  \text{separated  set for} \  T^{k_{\mu}} \right\} \leq \frac{\sharp \mathcal{F}_n}{\delta^2}$$

By Proposition \ref{main} (iii) it follows that :

\begin{eqnarray*}
& & h_{T^{k_{\mu}}}(M|F,\epsilon_{\mu}) \\
&\leq &  \lim_{n\rightarrow +\infty}\sup_{x\in 
F_{\sigma}} 
\left(2+H([\lambda^+_n(x,T^{k_{\mu}})-\chi^+(\nu,T^{k_{\mu}})]+3)\right)\left(\lambda^+_n(x,T^{k_{\mu}})-\chi^+(\nu,T^{k_{\mu}})\right)+ A
\end{eqnarray*}

According to the definition of $F_{\sigma}$  we have : 

\begin{eqnarray*}
\lim_{n\rightarrow +\infty}\inf_{x\in F_{\sigma}}\lambda^+_n(x,T^{k_{\mu}})-\chi^+(\nu,T^{k_{\mu}})&= & \lim_{n\rightarrow +\infty}\sup_{x\in F_{\sigma}}\lambda^+_n(x,T^{k_{\mu}})-\chi^+(\nu,T^{k_{\mu}})\\ 
& = & \int \log^+\|D_yT^{k_{\mu}}\|d\nu(y)- \chi^+(\nu,T^{k_{\mu}})\geq 0
\end{eqnarray*}

Moreover we deduce from the choice of $\delta_{\mu}$ that :

\begin{multline*}
\left|\int\log^+\|D_yT^{k_{\mu}}\|d\mu(y)- \chi^+(\nu,T^{k_{\mu}})\right|-\gamma  \leq \int \log^+\|D_yT^{k_{\mu}}\|d\nu(y)- \chi^+(\nu,T^{k_{\mu}})  \leq \\  \left|\int\log^+\|D_yT^{k_{\mu}}\|d\mu(y)- \chi^+(\nu,T^{k_{\mu}})\right|+\gamma
\end{multline*}

 and then by the choice of $k_{\mu}$ and since $\overline{\chi_0^+}(\mu,T)>\chi_0^+(\nu,T)-\gamma$ we get :

\begin{eqnarray*}
 \left|\overline{\chi^+_0}(\mu,T^{k_{\mu}})-\chi^+_0(\nu,T^{k_{\mu}})\right|-2k_{\mu}\gamma &\leq \int \log^+\|D_yT^{k_{\mu}}\|d\nu(y)- \chi^+_0(\nu,T^{k_{\mu}})& \leq \\
 &    \left|\overline{\chi^+_0}(\mu,T^{k_{\mu}})-\chi^+_0(\nu,T^{k_{\mu}})\right|+2k_{\mu}\gamma & \leq \\
 &   \overline{\chi^+_0}(\mu,T^{k_{\mu}})-\chi^+_0(\nu,T^{k_{\mu}})+4k_{\mu}\gamma
\end{eqnarray*}
 
Now we distinguish cases : 

\begin{itemize}
\item either $\left|\overline{\chi^+_0}(\mu,T)-\chi^+_0(\nu,T)\right|<3\gamma$, then  the term $\lim_{n\rightarrow +\infty}\sup_{x\in 
F_{\sigma}} H([\lambda^+_n(x,T^{k_{\mu}})-\chi^+_0(\nu,T^{k_{\mu}})]+3)\left(\lambda^+_n(x,T^{k_{\mu}})-\chi^+_0(\nu,T^{k_{\mu}})\right)$
is bounded above by $5\log(2)k_{\mu}\gamma$.

\item or $\overline{\chi^+_0}(\mu,T)\geq\chi^+_0(\nu,T)+3\gamma$, then we have $\lim_{n\rightarrow +\infty}\inf_{x\in F_{\sigma}}\lambda^+_n(x,T^{k_{\mu}})-\chi^+_0(\nu,T^{k_{\mu}})\geq k_{\mu}\gamma$. But recall we choose $k_{\mu}$ large enough so that $2R(T)H(k_{\mu}\gamma)<\gamma$. It follows that :
\begin{eqnarray*}
&\lim_{n\rightarrow +\infty} \sup_{x\in F_{\sigma}}H([\lambda^+_n(x,T^{k_{\mu}})-\chi^+_0(\nu,T^{k_{\mu}})]+3)\left(\lambda^+_n(x,T^{k_{\mu}})-\chi^+_0(\nu,T^{k_{\mu}})\right)  \\
\leq & H(k_{\mu}\gamma)2k_{\mu}R(T) \\ 
\leq & k_{\mu}\gamma 
\end{eqnarray*} 
\end{itemize}

We get finally in both cases  :

$$ h_{T^{k_{\mu}}}(M|F,\epsilon_{\mu})\leq  2\left( \overline{\chi_0^+}(\mu,T^{k_{\mu}})-\chi^+_0(\nu,T^{k_{\mu}})\right)+(5\log(2)+8)k_{\mu}\gamma+A $$

Then  by  letting $\sigma$ go to $1$ we obtain since $\frac{A}{k_{\mu}}<\gamma$ :

$$\frac{h_{T^{k_{\mu}}}^{New}(M|\nu,\epsilon_{\mu})}{k_{\mu}}\leq  
2\left(\overline{\chi_0^+}(\mu,T)-\chi^+_0(\nu,T)\right)+(5\log(2)+9)\gamma $$

Moreover according to Lemma \ref{tecnew}, we have $h_{T}^{New}(M|\nu,\epsilon_{\mu})\leq \frac{h_{T^{k_{\mu}}}^{New}(M|\nu,\epsilon_{\mu})}{k_{\mu}}$.
Therefore : 
$$ h_T^{New}(M|\nu,\epsilon_{\mu})\leq  
2\left(\overline{\chi_0^+}(\mu,T)-\chi^+_0(\nu,T)\right)+(5\log(2)+9)\gamma $$

This concludes the proof of Theorem \ref{super} because $\gamma$ can be chosen arbitrarily small.
\end{demof}

\section{Finite time stable fields}

We endow $\R^2$ with the Euclidian norm $\|\|$. The induced norm on the spaces $\mathcal{L}(\R^2,\R^2)$ and $\mathcal{L}(\mathcal{L}(\R^2,\R^2),\R^2)$ will be also denoted $\|\|$. Also we write $\|DT\|:=\sup_{x\in U}\|D_xT\|$ and $\|D_{\phi}T\|:=\sup_{x\in [0,1]^2}\|D_{\phi(x)}T\|$ where $\phi:[0,1]^2\rightarrow \R^2$ and $T$ is a $\mathcal{C}^1$ map defined on an open subset $U$ of $\R^2$ containing the image of $\phi$.

In the two following sections we consider a sequence $\mathcal{T}:=(T_n)_{n\in\N}$ of $\mathcal{C}^2$ diffeomorphisms from $B(0,2)\subset \R^2$ to $\R^2$ with  
$T_n(0)=0$ for all $n\in \N$.
For each $n\in \N$ we  denote $T^n$ the composition $T_n\circ...\circ T_1$ defined on $B(n,2):=\{x\in \R^2,\ 0\leq \forall k< n,\  \|T^kx\|<2\}$. 
 \\

We extend the notion of finite time hyperbolic set to this background by defining for all real numbers $\chi^+>0>\chi^-$, $\min(\chi^+,-\chi^-)>\gamma>0$ and $C>1$ the set :
\begin{multline*}
\mathcal{H}^n_{\mathcal{T}}(\chi^+,\chi^-,\gamma,C):=\{x\in B(n,2) \ : \ \forall 1\leq k\leq n, \ C^{-1}e^{(\chi^+-\gamma)k}\leq \|D_xT^k\|\leq Ce^{(\chi^++\gamma)k}, \\  
C^{-1}e^{(-\chi^--\gamma)k}\leq \|D_{T^kx}T^{-k}\|\leq Ce^{(-\chi^-+\gamma)k}\}
\end{multline*}

Recall a matrix $A\in \mathcal{G}l_2(\R)$ is said hyperbolic if $\|A\|,\|A^{-1}\|>1$. Such a matrix sends a circle to an ellipse. This defines two specific orthogonal  directions with orthogonal images : the most contracted and the most expanded one. Let $e_A$ and $f_A$ be unit vectors with arbitrary sense and with directions coinciding respectively with the most contracted and the most expanded directions, i.e. $\|e_A\|=\|f_A\|=1$, $\|Ae_A\|=\|A^{-1}\|^{-1}$ and $\|Af_A\|=\|A\|$. If $u,v\in \R^2$ we will denote by $\angle u,v$ the oriented angle of the two directed lines generated by $u$ and $v$, that is the angle belonging to $]-\pi,\pi]$  determined up to integral multiples of $2\pi$. We also write $|\angle|u,v:=|\angle u,v|$.\\

Let us denote by $U_n$ the open subset of $\R^2$ defined by $U_n:=\{x\in B(n,2), \ D_xT^n$ is hyperbolic$\}$, $e_n:U_n\rightarrow\R^2$ the most contracted field (also called the $n$ finite time stable field), that is $e_n(x)=e_{D_xT^n}$ for all $x\in U_n$, and $f_n:U_n\rightarrow \R^2$ 
the most expanded one (also called the $n$ finite time unstable field), that is $f_n(x)=f_{D_xT^n}$. One can choose the vector fields $e_n$ and $f_n$ to be continuous (an open subset of $\R^2$ is orientable) and $\angle e_n(x),f_n(x)=\frac{\pi}{2}$ for all $x\in U_n$. We also 
consider the vector fields $e_n^{(n)}(x)=\frac{D_xT^ne_n(x)}{\|D_xT^ne_n(x)\|}$ and 
$f_n^{(n)}(x)=\frac{D_xT^nf_n(x)}{\|D_xT^nf_n(x)\|}$ on $U_n$. Observe that $e_n^{(n)}\circ T^{-n}$ and $f_n^{(n)}\circ T^{-n}$ defined on $T^nU_n$ 
are respectively the expanded and the contracted fields for $T^{-n}$.  Since $T^n$ is a $\mathcal{C}^2$ map, all these vectors fields are in fact  
$\mathcal{C}^1$.

\scalebox{0.43}{\includegraphics{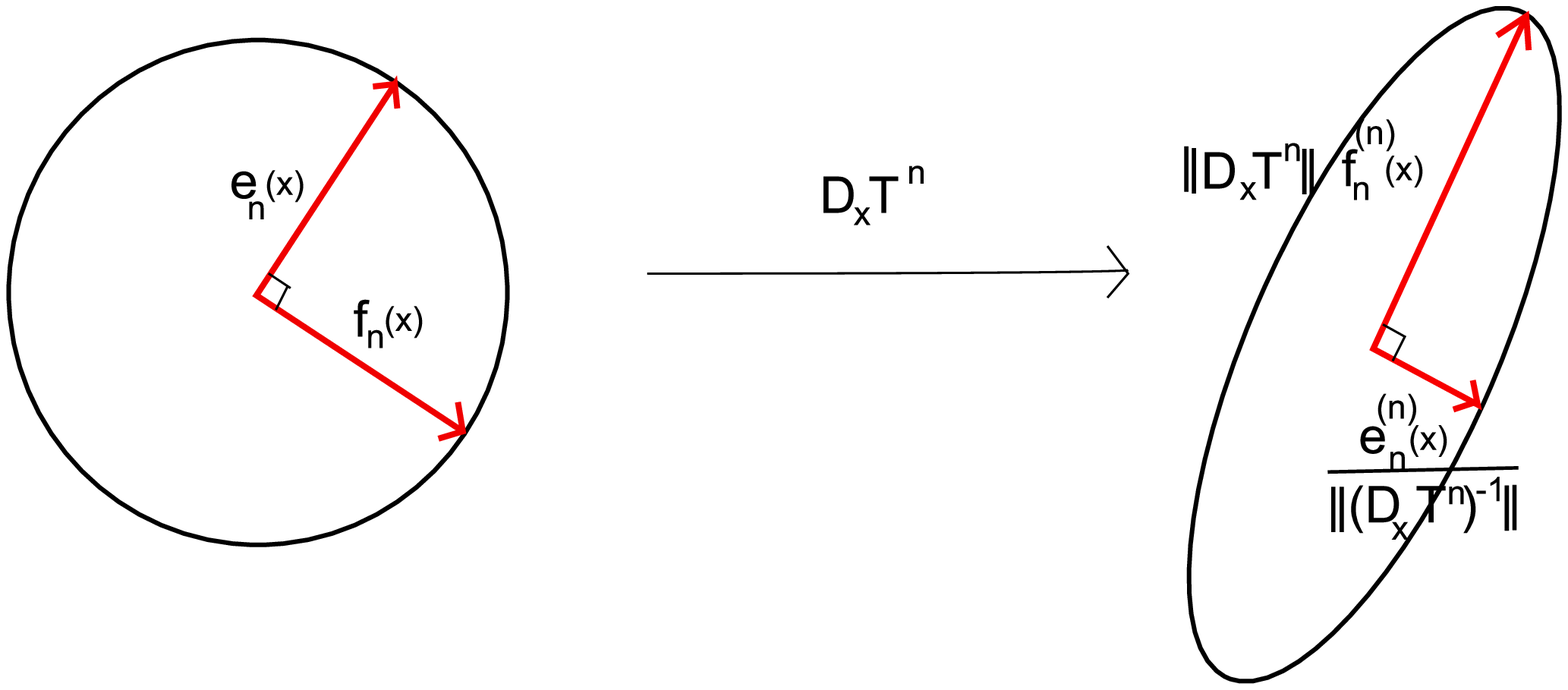}}

The vector fields $e_n$ and $f_n$ (resp. $e_n^{(n)}$ and $f_n^{(n)}$) are orthogonal so that, the angles $\angle e_n(x),e_n(y)$ (resp. $\angle 
e_n^{(n)}(x),e_n^{(n)}(y)$) and  $\angle f_n(x),f_n(y)$ (resp. $\angle f_n^{(n)}(x),f_n^{(n)}(y)$) are equal. As $e_n$ is a unit vector field, 
the image of $D_xe_n$ is orthogonal to $e_n(x)$ for all $x\in U_n$ (by derivating the relation $\|e_n(x)\|^2=1$). The same remark apply to 
$e_n^{(n)}$, $f_n$, $f_n^{(n)}$. Also by derivating the orthoganility relations $e_n(x).f_n(x)=0$ and $e_n^{(n)}(x).f_n^{(n)}(x)=0$ we get 
$D_xe_n(u).f_n(x)=-D_xf_n(u).e_n(x)$ and $D_xe_n^{(n)}(u).f_n=-D_xf_n^{(n)}(u).e_n^{(n)}(x)$. It follows that  
$\|D_xe_n(u)\|=\|D_xf_n(u)\|$ and $\|D_xe_n^{(n)}(u)\|=\|D_xf_n^{(n)}(u)\|$ for all $u\in \R^2$ and for all $x\in U_n$.\\

If $\phi:[0,1]^2\rightarrow \R^2$ is a $\mathcal{C}^1$ map, we say that $\phi$ satisfies $(H_n)$ when : 

\begin{eqnarray*} (\textbf{H}_\textbf{n}):   \ \ & \|D(x\mapsto D_{\phi(x)}T^n)\|  \leq K\|D_{\phi}T^n\| \ ; \\
& \|D(x\mapsto D_{T^n\circ\phi(x)}T^{-n})\| \leq  K\|D_{T^n\circ\phi}T^{-n}\|. 
\end{eqnarray*}

The real number $K$ is a very small universal constant ($K=10^{-10}$ is suitable). This constant will be slightly changed at finitely many steps of the proof but we always denote it by $K$ to simplify the computations. For example we will write $2K=K$, $e^{K}=1+K$, etc.   

Observe that $(H_n)$ is symmetrical by inversing the dynamical system, i.e. if $\phi:[0,1]^2\rightarrow \R^2$ is a $\mathcal{C}^1$ map satisfying $(H_n)$ for the map $T^n$, then so does $T^{n}\circ \phi$ for the map $T^{-n}$.\\

This assumption is the main key point of this paper. For interval maps it follows easily from the total order on $\R$ that the maximal cardinality of $(n,\delta)$ separated sets lying in a monotone branch of $f^n$ is less than $n/\delta$. Then to estimate the Newhouse local entropy one only needs to count the number of monotone branches (which coincide with the invertible branches) \cite{Dow}. Contrarily to the one dimensional case where positive entropy requires noninvertibility,  surface diffeomorphisms may have positive  entropy. The Smale horseshoe is a key model of such maps : indeed any surface diffeomorphism admits an invariant compact hyperbolic set in the closure of periodic points with topological entropy arbitrarily close to the global topological entropy \cite{Kat} and horseshoes then arise near homoclinic points according to Smale-Birkhoff homoclinic theorem \cite{guck}. In this context entropy  results in some sense from a combination of expansion, contraction and bending. One can try to define in dimension two a notion analagous to  the notion of monotone branches in dimension one. The natural idea is to consider pieces of the surface where the effect of bending is small and this corresponds exactly to the assumption $(H_n)$ (see the following proposition). For example the Smale horseshoe clearly does not satisfy $(H_n)$. Observe also that if $T$ is an interval map and $\phi:[0,1]\rightarrow [0,1]$ is a $\mathcal{C}^1$ map satisfying  $\|D(x\mapsto D_{\phi(x)}T^n)\|  \leq K\|D_{\phi}T^n\|$ then the image of $\phi$ lies in a monotone branch of $T^n$. 

 Another approach was developped by Y.Yomdin in the eighties \cite{Yoma} \cite{Yomb}. He used semi-algebraic tools to bound the local dynamical complexity of $\mathcal{C}^r$ maps. In Yomdin's theory which applies in any dimension and in any intermediate regularity one reparametrizes a Bowen ball of order $n$ by $C^r$ maps  $\phi$ defined on the unit square such that the $\mathcal{C}^r$ norm (i.e. the maximum over $1\leq s\leq r$ of the supremum norm of the $s^{th}$ derivative) of $T^k\circ \phi$ is less than $1$ for $0\leq k\leq n-1$ (in particular the maps $T^k\circ \phi$  are contracting and therefore  the reparametrization maps $\phi$ do not separate points). Remark that these conditions for $r\geq 2$ do not lead to the assumption $(H_n)$. Indeed  the reparametrization $\phi$ can be so contracting that the ocillation of $DT^n$ on the image of $\phi$ is still not small compared to its supremum norm $\|DT^n\|$. Conversely the assumption $(H_n)$ do not imply that $\|D^2(T\circ \phi)\|<1$. However to bound the entropy we do not need to control the second derivative of $T\circ \phi$ : as a naive example an invertible $\mathcal{C}^2$ interval map  do not separate points and may have a large $\mathcal{C}^2$ norm but the oscillation of the derivative is less than the supremum norm of the derivative to ensure the invertibility.\\

We introduce now a very weak assumption of hyperbolicity. If $\phi:[0,1]^2\rightarrow \R^2$ is a $\mathcal{C}^1$ map, we say that $\phi$ satisfies $(F_n)$ when  :
 
\begin{eqnarray*} (\textbf{F}_\textbf{n}): & \forall t\in [0,1]^2,  \ \ \|D_{\phi(t)}T^n\|\|D_{T^n\circ\phi(t)}T^{-n}\|\geq 2 \ ; \\ 
 & \|D_{\phi}T^n\|\|D_{T^n\circ\phi}T^{-n}\|\geq \max_{k=0,...,n}\|D_{\phi}T^k\|.   
\end{eqnarray*} 

In the following proposition we prove that under the assumptions $(H_n)$ and $(F_n)$ the finite time stable and unstable fields $e_n,f_n$ and their image $e_n^{(n)},f_n^{(n)}$ do not oscillate on the image of $\phi$. 

\begin{prop}\label{angles}
Let $n\in\N$ and let $\phi:[0,1]^2\rightarrow \R^2$ be a $\mathcal{C}^1$ map 
satisfying  $(H_n)$ and $(F_n)$ then 
 we have  for all $y\in [0,1]^2$ : $$\|D_y(e_n\circ \phi)\|\leq K$$
$$\|D_y(e_n^{(n)}\circ \phi)\|\leq K$$
\end{prop}

\begin{demo}
Clearly it is enough to prove that for all $t\in [0,1]^2$ there exists a function $\epsilon:\R^+\rightarrow \R^+$ with $\lim_{t\rightarrow 0}\epsilon(t)=0$ such that $|\angle| e_n\circ\phi(t),e_n\circ \phi(s)\leq K\|s-t\| +\|s-t\|\epsilon(\|s-t\|)$ and 
$|\angle| e_n^{(n)}\circ\phi(t),e_n^{(n)}\circ \phi(s)\leq K\|s-t\| +\|s-t\|\epsilon(\|s-t\|)$ for all $s \in [0,1]^2$, that is with the Laudau notation $|\angle| e_n\circ\phi(t),e_n\circ \phi(s)\leq K\|s-t\| +o_{s\rightarrow t}(\|s-t\|)$ and
$|\angle| e_n^{(n)}\circ\phi(t),e_n^{(n)}\circ \phi(s)\leq K\|s-t\| +o_{s\rightarrow t}(\|s-t\|)$.\\

Let  $t,s\in [0,1]^2$. It follows from the triangular inequality that :

\begin{eqnarray*}
|\angle| e_n^{(n)}\circ\phi(t),e_n^{(n)}\circ \phi(s)& =& |\angle| f_n^{(n)}\circ\phi(t),f_n^{(n)}\circ \phi(s)\\
&= &|\angle| D_{\phi(t)}T^nf_n(\phi(t)),D_{\phi(s)}T^nf_n(\phi(s))\\
&\leq &|\angle| D_{\phi(t)}T^nf_n(\phi(t)),D_{\phi(t)}T^nf_n(\phi(s)) \\
& &\ \ \ \  \ \ \ \ \  \ \ \ \ \ \ \ \ \ \ \ \ \ + |\angle|  D_{\phi(t)}T^nf_n(\phi(s)),D_{\phi(s)}T^nf_n(\phi(s))
\end{eqnarray*}

The first part of the right member can be bounded in the following way : 
$$\tan  \angle D_{\phi(t)}T^nf_n(\phi(t)),D_{\phi(t)}T^nf_n(\phi(s)) = \frac{\tan\angle f_n(\phi(t)),f_n(\phi(s))}{\|D_{\phi(t)}T^n\|\|D_{T^n\phi(t)}T^{-n}\|}$$

Indeed we have $f_n(\phi(s))=f_n(\phi(t))\cos \angle f_n(\phi(s)),f_n(\phi(t))+e_n(\phi(t))\sin \angle f_n(\phi(s)),f_n(\phi(t))$. 
Then by applying $D_{\phi(t)}T^n$ :
\begin{multline*}
D_{\phi(t)}T^nf_n(\phi(s))=f_n^{(n)}(\phi(t))\|D_{\phi(t)}T^n\|\cos\angle f_n(\phi(t)),f_n(\phi(s))+\\e^{(n)}_n(\phi(t))\|D_{T^n\circ \phi(t)}T^{-n}\|^{-1}\sin\angle f_n(\phi(t)),f_n(\phi(s))
\end{multline*}

 Therefore
 
\begin{eqnarray*}
\tan  \angle D_{\phi(t)}T^nf_n(\phi(t)),D_{\phi(t)}T^nf_n(\phi(s)) &=& \frac{\|D_{T^n\circ \phi(t)}T^{-n}\|^{-1}\sin\angle f_n(\phi(t)),f_n(\phi(s))}{\|D_{\phi(t)}T^n\|\cos\angle f_n(\phi(t)),f_n(\phi(s))}\\
&=&\frac{\tan\angle f_n(\phi(t)),f_n(\phi(s))}{\|D_{\phi(t)}T^n\|\|D_{T^n\circ \phi(t)}T^{-n}\|}
\end{eqnarray*}

Consider now the second part $|\angle|  D_{\phi(t)}T^nf_n(\phi(s)),D_{\phi(s)}T^nf_n(\phi(s))$. We have :

\begin{eqnarray*}
\|D_{\phi(t)}T^nf_n(\phi(s))-D_{\phi(s)}T^nf_n(\phi(s))\| & \geq & |D_{\phi(t)}T^nf_n(\phi(s)).e_n^{(n)}(\phi(s))|\\
& \geq & \|D_{\phi(t)}T^nf_n(\phi(s))\| \times \\
& &   \  \ \  \ \ \ \ \sin|\angle| D_{\phi(t)}T^nf_n({\phi(s)}),D_{\phi(s)}T^nf_n(\phi(s))
\end{eqnarray*}

But we have $K\|s-t\|\|D_{\phi}T^n\|\geq \|D_{\phi(t)}T^nf_n(\phi(s))-D_{\phi(s)}T^nf_n(\phi(s))\|$ and $\|D_{\phi(t)}T^n\|\geq (1-K)\|D_{\phi}T^n\|$  by assumption ($H_n$). Moreover by continuity $\|D_{\phi(t)}T^nf_n(\phi(s))\|$ goes to $\|D_{\phi(t)}T^n\|$ when $s$ goes to $t$. It implies that $\sin |\angle| D_{\phi(t)}T^nf_n(\phi(s)), D_{\phi(s)}T^nf_n(\phi(s))\leq K\|s-t\|+o_{s\rightarrow t}(\|s-t\|)$.
Finally we have : 
$$|\angle| e_n^{(n)}\circ\phi(t),e_n^{(n)}\circ \phi(s)\leq \frac{|\angle| f_n(\phi(t)),f_n(\phi(s))}{\|D_{\phi(t)}T^n\|\|D_{T^n\circ\phi(t)}T^{-n}\|} + K\|s-t\|+o_{s\rightarrow t}(\|s-t\|)$$

By considering $T^n\circ\phi$ instead of $\phi$ and $T^{-1}$ instead of $T$ we get symetrically : 
 
$$|\angle| f_n\circ\phi(t),f_n\circ \phi(s)\leq \frac{|\angle| e_n^{(n)}(\phi(t)),e_n^{(n)}(\phi(s))}{\|D_{\phi(t)}T^n\|\|D_{T^n\circ\phi(t)}T^{-n}\|} + K\|s-t\|+o_{s\rightarrow t}(\|s-t\|)$$

But by assumption $(F_n)$ we have $\|D_{\phi(t)}T^n\|\|D_{T^n\circ\phi(t)}T^{-n}\|\geq 2$ for all $t\in [0,1]^2$ and then one concludes the proof of the lemma by combining the  two previous inequalities.

\end{demo}

Observe that in fact we only need to suppose the first part of the assumption $(F_n)$ (which is symmetrical) in the previous proof ; the second part will be useful in the next sections.

\section{Finite time Rectangle}

 The vectors fields $e_n,f_n:U_n\rightarrow \R^2$ are $\mathcal{C}^1$ so that they can be locally 
 integrated. We consider the integral curves $\mathcal{E}_n(z)$ (resp. $\mathcal{F}_n(z)$) through $z$ to the field $e_n$ (resp. $f_n$). These 
 curves are called $n$ finite time stable (resp. unstable) manifolds. This approach comes from the works of M.Benedicks and L.Carleson on the Henon map \cite{ben} and was formalized by M.Holland and S.Luzzato \cite{hol} \cite{holl}. In particular, under some assumptions of hyperbolicity on surfaces, M.Holland and S.Luzzato \cite{hol} give an alternative proof of the existence of stable manifolds by proving the convergence of finite time stable manifolds.  Following the notion of rectangle in the theory of hyperbolic dynamical systems we introduce now finite time rectangles. This new finite time concept is particularly well adapted to the computation of entropy. \\

\begin{defi}
A subset $R_n$ of  $U_n$ is called a $n$-rectangle if there exists a diffeomorphism  $\phi:V\rightarrow W\subset U_n$ defined on an open neighborhood  $V$ of $[0,1]^2$ with $\phi([0,1]^2)=R_n$, such that $\phi([0,1]\times \{0\})$ and $\phi([0,1] \times  \{1\})$ (resp.  $\phi(\{0\}\times [0,1])$ and $\phi(\{1\}\times [0,1])$ are included in a $n$ finite time stable (resp. unstable) manifold. 
\end{defi}

In the following the map $\phi$ and the couple $(\phi,R_n)$ will also refer to a $n$-rectangle.
Observe that if $(\phi,R_n)$ is a $n$-rectangle for $T^n$ then $(T^n\circ \phi,T^{n}R_n)$ is a $n$-rectangle for $T^{-n}$.
\\

We can give a satisfactory geometrical reparametrization of a $n$-rectangle $R_n$. For any point $x$ of $R_n$, we reparametrize $\mathcal{E}_n(x)\cap R_n$ and $\mathcal{F}_n(x)\cap R_n$ from $[0,1]$ with constant rate in the same direction as $e_n(x)$ and $f_n(x)$. We get in this way two maps of the interval $\phi^e_x:[0,1]\rightarrow \mathcal{E}_n(x)\cap R_n$ and $\phi^f_x:[0,1]\rightarrow \mathcal{F}_n(x)\cap R_n$. Then the map $\phi_{x,n}:[0,1]^2\rightarrow R_n$ defined by

$$\phi_{x,n}(t,s)=\mathcal{F}_n(\phi_x^e(t))\cap \mathcal{E}_n(\phi_x^f(s))$$
is a diffeomorphism from $[0,1]^2$ to $M$. One can  associate to any subrectangle $[a,b]\times [c,d]$ of the unit square the $n$-rectangle $\phi_{x,n}([a,b]\times[c,d])$. Remark also that 
$\partial_1\phi_{x,n}(t,s)$ (resp. $\partial_2\phi_{x,n}(t,s)$) is colinear to $e_n(\phi_{x,n}(t,s))$ (resp. to $f_n(\phi_{x,n}(t,s))$). In fact $\phi_{x,n}$ is a foliation box for the foliations $\mathcal{F}_n$ and $\mathcal{E}_n$ simultaneously. Such reparametrizations $\phi_{x,n}$ will be called admissible charts.\\

 To simplify the notations, we write $a\simeq^{\lambda} b$ for $a,b>0$ and $\lambda\geq1$  which means that  $\lambda^{-1}\leq \frac{a}{b}\leq \lambda$. Moreover if $\gamma$ is a $\mathcal{C}^1$ curve its lenght will be denoted by $l(\gamma)$.
We say that a $n$-rectangle $R_n$ satisfies  $(G_n)$ when :

\begin{eqnarray*} (\textbf{G}_\textbf{n}): & \forall x\in R_n,    & (1+K)l(\mathcal{E}_n(x)\cap R_n)\geq l(\mathcal{F}_n(x)\cap R_n)\geq (1-K)\frac{l(\mathcal{E}_n(x)\cap R_n)}{\max_{k=0,...,n}\|D_xT^k\|}.  
\end{eqnarray*}

The assumption $(G_n)$ is not symmetrical :  if $R_n$ is a $n$-rectangle satisfying $(G_n)$ then this assumption is not satisfied in general by the $n$-rectangle $T^nR_n$ for $T^{-n}$. However we prove in the next lemma that $(1+K)l(T^n\mathcal{F}_{n}(x)\cap T^nR_n)\geq l(T^n\mathcal{E}_n(x)\cap T^nR_n)$ assuming $(H_n)$, $(G_n)$ and  $(F_n)$ for $R_n$.

\subsection{Finite time rectangle with small oscillation of the derivative}

We show that under the assumptions $(H_n)$, $(G_n)$ and  $(F_n)$ the $n$-rectangles really look like rectangles : the $n$ finite time stable (resp. unstable) manifolds foliating a $n$-rectangle have almost the same length (Proposition \ref{marre}) and more precisely the map $\phi_{x,n}$ reparametrizes these manifolds almost with constant rate (Proposition \ref{angle}).   

\begin{lem}\label{zut}
Let $n\in\N$ and let $(\phi,R_n)$ be a $n$-rectangle satisfying $(H_n)$, $(G_n)$ and $(F_n)$. Then we have for all $x\in R_n$ :

$$(1+K)l(T^n\mathcal{F}_{n}(x)\cap T^nR_n)\geq l(T^n\mathcal{E}_n(x)\cap T^nR_n)$$
\end{lem}

\begin{demo}
 First remark that the assumption $(H_n)$, which states that $\|D(x\mapsto D_{\phi(x)}T^n)\|\leq K\|D_{\phi}T^n\|$ and $\|D(x\mapsto D_{T^n\circ\phi(x)}T^{-n})\|\leq K\|D_{T^n\circ\phi}T^{-n}\|$, imply respectively that $\|D_{x}T^n\|\simeq^{1+K}\|D_{\phi}T^n\|$ and $\|D_{T^n(x)}T^{-n}\|\simeq^{1+K}\|D_{T^n\circ\phi}T^{-n}\|$ for all $x\in R_n$. Therefore 

$$ l(T^n\mathcal{F}_{n}(x)\cap T^nR_n)\simeq^{1+K}\|D_{\phi}T^n\|l(\mathcal{F}_n(x)\cap R_n)$$ 
$$l(T^n\mathcal{E}_{n}(x)\cap T^nR_n)\simeq^{1+K}\|D_{T^n\circ\phi}T^{-n}\|^{-1}l(\mathcal{E}_n(x)\cap R_n)$$

and then according to the hypotheses $(G_n)$ and $(F_n)$ 
\begin{eqnarray*}
\frac{l(T^n\mathcal{F}_{n}(x)\cap T^nR_n)}{l(T^n\mathcal{E}_{n}(x)\cap T^nR_n)}&\simeq^{1+K}&\frac{\|D_{\phi}T^n\|\|D_{T^n\circ\phi}T^{-n}\|l(\mathcal{F}_n(x)\cap R_n)}{l(\mathcal{E}_n(x)\cap R_n)}\\
& \geq &\frac{\|D_{\phi}T^n\|\|D_{T^n\circ\phi}T^{-n}\|}{(1+K)\max_{k=0,...,n}\|D_xT^k\|}\\
&\geq & \frac{1}{1+K}
\end{eqnarray*}

%But since $R_n\subset \mathcal{H}_{\mathcal{T}}^n(\chi^+,\chi^-,\gamma,C)$, we have for large $n$ (recall $-\chi^->3\gamma$): 

%$$\frac{\|D_{\phi}T^n\|}{\|D_{T^n\circ\phi}T^{-n}\|\max_{k=0,...,n}\|D_xT^k\|}\geq C^{-3}e^{n(\chi^--3\gamma)}\geq \frac{1}{1-K}$$

This concludes the proof of the lemma.
\end{demo}

\begin{prop}\label{marre}
Let $(\phi,R_n)$ be a $n$-rectangle satisfying $(H_n)$, $(G_n)$ and $(F_n)$. Then we have for all $x,y\in R_n$ :

$$l(\mathcal{E}_n(x)\cap R_n)\simeq^{1+K}l(\mathcal{E}_n(y)\cap R_n)$$
$$l(\mathcal{F}_n(x)\cap R_n)\simeq^{1+K}l(\mathcal{F}_n(y)\cap R_n)$$

\end{prop}

\begin{demo}
Fix $x,y\in R_n$ and let $z\in R_n$ be the intersection point of $\mathcal{E}_n(y)$ and $\mathcal{F}_n(x)$. One deduces easily from $|\angle| e_n(u),e_n(v)\leq K$ for all $u,v\in R_n$ that :

\begin{center}
\begin{itemize}
\item $l(\mathcal{E}_n(x)\cap R_n)\simeq^{1+K}\|\phi^e_x(1)-x\|+\|x-\phi^e_x(0)\|$ ;
\item $l(\mathcal{E}_n(y)\cap R_n)\simeq^{1+K}\|\phi^e_y(1)-z\|+\|z-\phi^e_y(0)\|$.
\end{itemize}
\end{center}

Moreover the quadrilaterals with vertices $\phi^e_x(0),\phi^e_y(0),x,z$ and $\phi^e_x(1),\phi^e_y(1),x,z$  are almost rectangles : their interior angles differ from $\frac{\pi}{2}$ by at most $K$. Then  we have by trivial trigonometric arguments $\|x-\phi_x^e(1)\|\leq (1+K)\|z-\phi^e_y(1)\|+K\|x-z\|$ and $\|x-\phi_x^e(0)\|\leq (1+K)\|z-\phi^e_y(0)\|+K\|x-z\|$. Going back to the length of the finite time stable manifolds we obtain  $l(\mathcal{E}_n(x)\cap R_n)\leq (1+K)l(\mathcal{E}_n(y)\cap R_n)+Kl(\mathcal{F}_n(x)\cap R_n)$ and finally $l(\mathcal{E}_n(x)\cap R_n)\leq (1+K)l(\mathcal{E}_n(y)\cap R_n)$ since we have $(1+K)l(\mathcal{E}_n(x)\cap R_n)\geq l(\mathcal{F}_n(x)\cap R_n)$  by assumption $(F_n)$. This concludes the first point of the lemma  by switching the role of $x$ and $y$.

 One can apply the previous proof to the $n$-rectangle $T^nR_n$ for $T^{-1}$ because  $(1+K)l(T^n\mathcal{F}_{n}(x)\cap T^nR_n)\geq l(T^n\mathcal{E}_n(x)\cap T^nR_n)$ for all $x\in R_n$ according to the previous lemma.  We get then $l(T^n\mathcal{F}_{n}(x)\cap T^nR_n)\simeq^{1+K} 
l(T^n\mathcal{F}_{n}(y)\cap T^nR_n)$ and finally $l(\mathcal{F}_{n}(x)\cap R_n)\simeq^{1+K} 
l(\mathcal{F}_{n}(y)\cap R_n)$ since $l(T^n\mathcal{F}_{n}(u)\cap T^nR_n)\simeq^{1+K}\|D_{\phi}T^n\|l(\mathcal{F}_{n}(u)\cap R_n)$ for all $u\in R_n$ by assumption $(H_n)$. 

\scalebox{0.45}{\includegraphics{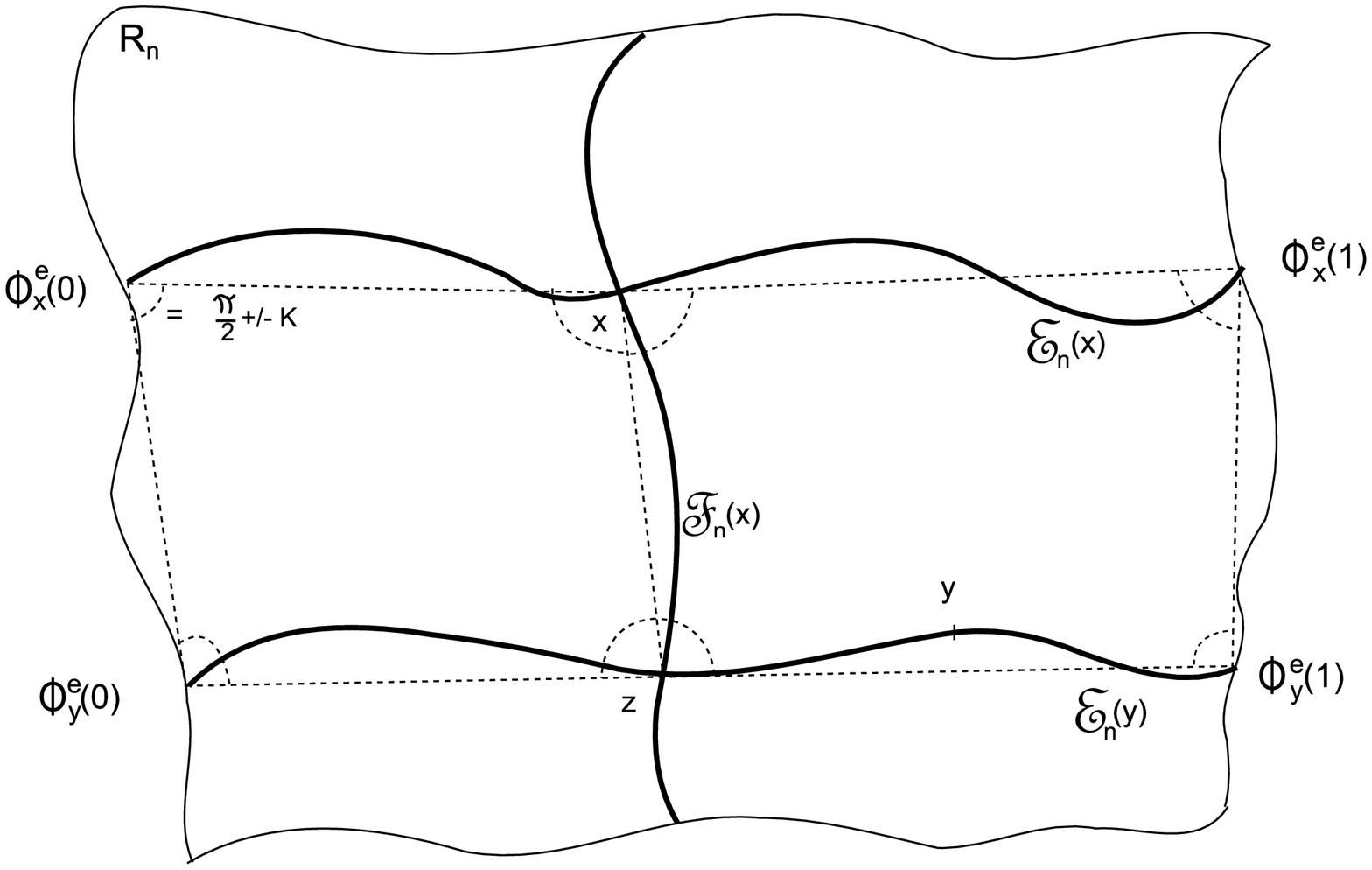}}

\end{demo}

In Lemma \ref{zut} and Proposition \ref{marre} we did not require that the reparametrization $\phi$ of the $n$-rectangle are admissible charts. From now we will only consider reparametrizations of $n$-rectangles  of this kind.

\begin{prop}\label{angle}
Let $(\phi_{x,n},R_n)$ be a $n$-rectangle satisfying $(H_n)$, $(F_n)$ and $(G_n)$. Then we have for all $t,s\in [0,1]^2$ : $$\|\partial_1\phi_{x,n}(t,s)\| \simeq^{1+K} l(\mathcal{E}_n(x)\cap R_n)$$
$$\|\partial_2\phi_{x,n}(t,s)\| \simeq^{1+K} l(\mathcal{F}_n(x)\cap R_n)$$
\end{prop}

\begin{demo}
We first prove that $\|\partial_1\phi_{x,n}(t,s)\| \simeq^{1+K} l(\mathcal{E}_n(x)\cap R_n)$. For simplicity we put $\phi:=\phi_{x,n}$. Let $s_0\in [0,1]$ such that $x\in \phi([0,1]\times s_0)$. Since the $n$ finite time stable manifold of $x$ is reparametrized with constant rate we have $\|\partial_1\phi(t,s_0)\|=\|\partial_1\phi^e_x(t)\|=l(\mathcal{E}_n(x)\cap R_n)$ for all $t\in [0,1]$. Fix $t,s,s'\in [0,1]$, we have :

\begin{eqnarray*}
\|\partial_1\phi(t,s)\|-\|\partial_1\phi(t,s')\| & = & \partial_1\phi(t,s).e_n(\phi(t,s))-\partial_1\phi(t,s').e_n(\phi(t,s')) \\
& = & \int_{s'}^{s} \partial_2\left(\partial_1\phi(t,u).e_n(\phi(t,u))\right)du\\
& = & \int_{s'}^{s} \partial_2\partial_1\phi(t,u).e_n(\phi(t,u))du+\int_{s'}^{s}\partial_1\phi(t,s).\partial_2 e_n(\phi(t,u))du 
\end{eqnarray*}

Now $\partial_2\partial_1\phi(t,u)=\partial_1\partial_2\phi(t,u)=\partial_1\|\partial_2\phi(t,u)\|f_n(\phi(t,u))+\|\partial_2\phi(t,u)
\|\partial_1f_n(\phi(t,u))$ so that by Proposition \ref{angles} :

\begin{eqnarray*} 
\left|\int_{s'}^{s} \partial_2\partial_1\phi(t,u).e_n(\phi(t,u))du\right|&\leq & K\int_{s'}^{s} \|\partial_2\phi(t,u)\|du\leq Kl(\mathcal{F}_n(\phi(t,s_0))\cap R_n)
\end{eqnarray*}

and then according to Proposition 3 and assumption $(G_n)$ :

\begin{eqnarray*}
\left|\int_{s'}^{s} \partial_2\partial_1\phi(t,u).e_n(\phi(t,u))du\right|&\leq &  Kl(\mathcal{F}_n(x)\cap R_n) \\
&\leq & Kl(\mathcal{E}_n(x)\cap R_n) 
\end{eqnarray*}

Moreover 

\begin{eqnarray*}
\left|\int_{s'}^{s}\partial_1\phi(t,u).\partial_2 e_n(\phi(t,u))du\right|&\leq & K\int_{s'}^{s}\|\partial_1\phi(t,u)\|du
\end{eqnarray*}

We deduce from the two previous inequalities that  : 

$$\|\partial_1\phi(t,s)\|\leq \|\partial_1\phi(t,s')\|+Kl(\mathcal{E}_n(x)\cap R_n)+K\int_{s'}^{s}\|\partial_1\phi(t,u)\|du$$

so that by Grönwall Lemma we get $||\partial_1\phi(t,s)\| \leq \left(\|\partial_1\phi(t,s')\|+Kl(\mathcal{E}_n(x)\cap R_n)\right) e^{K|s-s'|}\leq (1+K)\|\partial_1\phi(t,s')\|+Kl(\mathcal{E}_n(x)\cap R_n)$.  Then by choosing $s=s_0$ or $s'=s_0$ we conclude that 
 $\|\partial_1\phi(t,s)\|\simeq^{1+K}l(\mathcal{E}_n(x)\cap R_n)$ for all $(t,s) \in[0,1]^2$. To get the other relation $\|\partial_2\phi_{x,n}(t,s)\| \simeq^{1+K} l(\mathcal{F}_n(x)\cap R_n)$ we apply the previous proof to the $n$-rectangle $T^nR_n$ for $T^{-1}$.

\end{demo}

We deduce from the two previous propositions that the change of admissible charts of a given $n$-rectangle is $\mathcal{C}^1$ bounded so that the choice of the admissible chart do not interfere to bound the dynamical complexity of $T$ on $R_n$ : 

\begin{coro}
Let $R_n$ a $n$-rectangle satisfying $(H_n)$, $(G_n)$ and $(F_n)$.  Then we have for all $x,y\in R_n$ : 
$$\|D(\phi_{x,n}\circ \phi_{y,n}^{-1})\|\leq 1+K$$
\end{coro}

The following Corollary allows us to bound the derivative of $T^n\circ \phi_{x,n}$ as soon as the diameter of $T^nR_n$ is bounded, i.e. we control the derivative of this map by the size of its image. 

\begin{coro}\label{croco}
Let $(\phi_{x,n},R_n)$ be a $n$-rectangle satisfying $(H_n)$, $(G_n)$, $(F_n)$ and such that the diameter of $T^nR_n$ is less than $2$. Then we have for all $x\in R_n$ :

$$\|D(T^n\circ \phi_{x,n})\|\leq 2+K$$ 

\end{coro}

\begin{demo}
It's enough to bound the derivative in the unstable direction, that is $\|\partial_2 (T^n\circ \phi_{x,n})\|\leq 2+K$. 
We have for all $t\in [0,1]$ : $$\left\|\int_0^1\partial_2 (T^n\circ \phi_{x,n})(t,s)ds\right\| =  \|T^n\circ\phi_{x,n}(t,1)-T^n\circ\phi_{x,n}(t,0)\|\leq 2$$
Remark now that :

\begin{eqnarray*}
 (1-K)\int_0^1\|\partial_2 (T^n\circ \phi_{x,n})(t,s)\|ds \leq & \int_0^1 \partial_2 (T^n\circ \phi_{x,n})(t,s).f_n^{(n)}(x)ds\leq &\left\|\int_0^1\partial_2 (T^n\circ \phi_{x,n})(t,s)ds\right\|
\end{eqnarray*}

We have also $\|\partial_2 (T^n\circ \phi_{x,n})(t,s)\|=\|\partial_2\phi_{x,n}(t,s)\|\times \|D_{\phi_{x,n}(t,s)}T^n\|\simeq^{1+K} 
l(\mathcal{F}_n(x)\cap R_n)\|D_xT^n\|$ for all $(t,s)\in [0,1]^2$ by Proposition \ref{angle}. Therefore $\|\partial_2 (T^n\circ \phi_{x,n})(t,s)\|\leq 2+K$ for all $(t,s) \in [0,1]^2$. 

\end{demo}

\subsection{From $n$- to $n+1$-rectangles}\label{supsub}

The end of this section devoted to the analysis of finite time rectangles deals with the links between $n$-rectangles and $n+1$-rectangles. This point is crucial for the induction step in the proof of Proposition \ref{main} presented in the next section.\\
 
 In the following we consider  sequences of maps $\mathcal{T}=(T_n)_{n\in\N}$ such that $T_n$ and its inverse $T_n^{-1}$ is uniformly $\mathcal{C}^1$ bounded, that is $S(\mathcal{T}):=\max(\sup_n\|DT_n\|, \sup_{n}\|DT_n^{-1}\|)<+\infty$. 

To compare $n$- and $n+1$-rectangles we also need stronger hyperbolicity assumptions than the weak hypothesis $(F_n)$. Let $\chi^+>0>\chi^-, \frac{\min(\chi^+,-\chi^-)}{3}>\gamma>0$ and $C>1$. We will consider $n$-rectangles $(\phi,R_n)$ included  in the finite time hyperbolic set $\mathcal{H}^n_{\mathcal{T}}(\chi^+,\chi^-,\gamma,C)$ with "$n$ large". By $n$ large we mean the statements of this subsection hold for $n$ larger than an integer $N$ depending only on $\chi^+,\chi^-,\gamma,C$ and $S(\mathcal{T})$. 

 Oberve that such $n$-rectangles satisfy the assumption $(F_n)$ and even $(F_{n+1})$. Indeed if  $(\phi,R_n)$ is a $n$-rectangle with $R_n\subset \mathcal{H}_{\mathcal{T}}^n(\chi^+,\chi^-,\gamma,C)$, then we have for large $n$ : 

$$\frac{\|D_{\phi}T^n\|}{\|D_{T^n\circ\phi}T^{-n}\|\max_{k=0,...,n}\|D_xT^k\|}\geq C^{-3}e^{n(-\chi^--3\gamma)}\geq 1$$

and for all $t\in [0,1]^2$ :

$$\|D_{\phi(t)}T^n\|\|D_{T^n\circ\phi(t)}T^{-n}\| \geq C^{-2}e^{n(\chi^+-\chi^--2\gamma)}\geq 2$$

The  next lemma estimates the angle $|\angle| e_n(x),e_{n+1}(x)$ under the condition of finite time hyperbolicity :

\begin{lem}\label{holangle}
For all $x\in \mathcal{H}^n_{\mathcal{T}}(\chi^+,\chi^-,\gamma,C)$ and for large $n$, 

\begin{eqnarray}
\tan|\angle | e_n(x),e_{n+1}(x)\leq \frac{2\|D_{T^nx}T_{n+1}\| \|(D_{T^nx}T_{n+1})^{-1}\|}{\|(D_xT^n)^{-1}\|\|D_xT^n\|}
\end{eqnarray}

\end{lem}

\begin{demo}
We write $e_{n}(x)=e_{n+1}(x)\cos \angle e_n(x),e_{n+1}(x) +f_{n+1}(x)\sin\angle e_n(x),e_{n+1}(x)$ and we apply $D_xT^{n+1}$ :

\begin{multline*}
D_xT^{n+1}e_n(x)=e_{n+1}^{(n+1)}(x)\|D_{T^{n+1}x}T^{-n-1}\|^{-1}\cos \angle e_n(x),e_{n+1}(x)+\\
f_{n+1}^{(n+1)}(x)\|D_xT^{n+1}\|\sin\angle e_n(x),e_{n+1}(x)
\end{multline*}

and then by taking the Euclidian norm :
$$\|D_xT^{n+1}e_n(x)\|^2=\|D_{T^{n+1}x}T^{-n-1}\|^{-2}\cos^2 \angle e_n(x),e_{n+1}(x)+\|D_xT^{n+1}\|^2\sin^2\angle e_n(x),e_{n+1}(x)$$

We compute :

\begin{eqnarray*}
\tan|\angle| e_n(x),e_{n+1}(x) &=&\left(\frac{\|D_xT^{n+1}e_n(x)\|^2-\|D_{T^{n+1}x}T^{-n-1}\|^{-2}}{\|D_xT^{n+1}\|^2-\|D_xT^{n+1}e_n(x)\|^2}\right)^{\frac{1}{2}}
\end{eqnarray*}

 As $x\in \mathcal{H}^n_{\mathcal{T}}(\chi^+,\chi^-,\gamma,C)$ we get for $n$ large enough :
\begin{eqnarray*}
\tan|\angle| e_n(x),e_{n+1}(x)&\leq &\frac{\|D_xT^{n+1}e_n(x)\|/\|D_xT^{n+1}\|}{\left(1-\|D_xT^{n+1}e_n(x)\|^2/\|D_xT^{n+1}\|^2\right)^{\frac{1}{2}}}\\
&\leq & 2\frac{\|D_xT^{n+1}e_n(x)\|}{\|D_xT^{n+1}\|}\\
&\leq & \frac{2\|D_{T^nx}T_{n+1}\|\|(D_{T^nx}T_{n+1})^{-1}\|}{\|(D_xT^n)^{-1}\|\|D_xT^n\|}
\end{eqnarray*}

This concludes the proof of the lemma. 
\end{demo}

Given a $n$-rectangle $(\phi_{x,n},R_n)$ we define its saturated set, $Sat(R_n)$, as the smallest subset of $R_n$ saturated  in $n+1$ stable and unstable manifolds containing the middle third of $R_n$, that is \footnote{Note that $Sat(R_n)$ depends on the choice of the admissible chart $\phi_{x,n}$} :

$$Sat(R_n):=\left\{z\in R_n, \ \exists u,v\in \phi_{x,n}\left(\left[\frac{1}{3},\frac{2}{3}\right]^2\right), \ s.t.\ z\in \mathcal{E}_{n+1}(u)\cap \mathcal{F}_{n+1}(v) \right\}$$

\begin{prop}\label{sat}
For any $n$-rectangle $(\phi_{x,n},R_n)$ with large $n$ included in 
$\mathcal{H}_{\mathcal{T}}^n(\chi^+,\chi^-,\gamma,C)$ satisfying $(H_n)$ and $(G_n)$, the saturated set $Sat(R_n)$ of $R_n$ 
 is a $n+1$-rectangle 
included in $\phi_{x,n}\left(\left[\frac{1}{3}-K,\frac{2}{3}+K\right]^2\right)$.

\scalebox{0.43}{\includegraphics{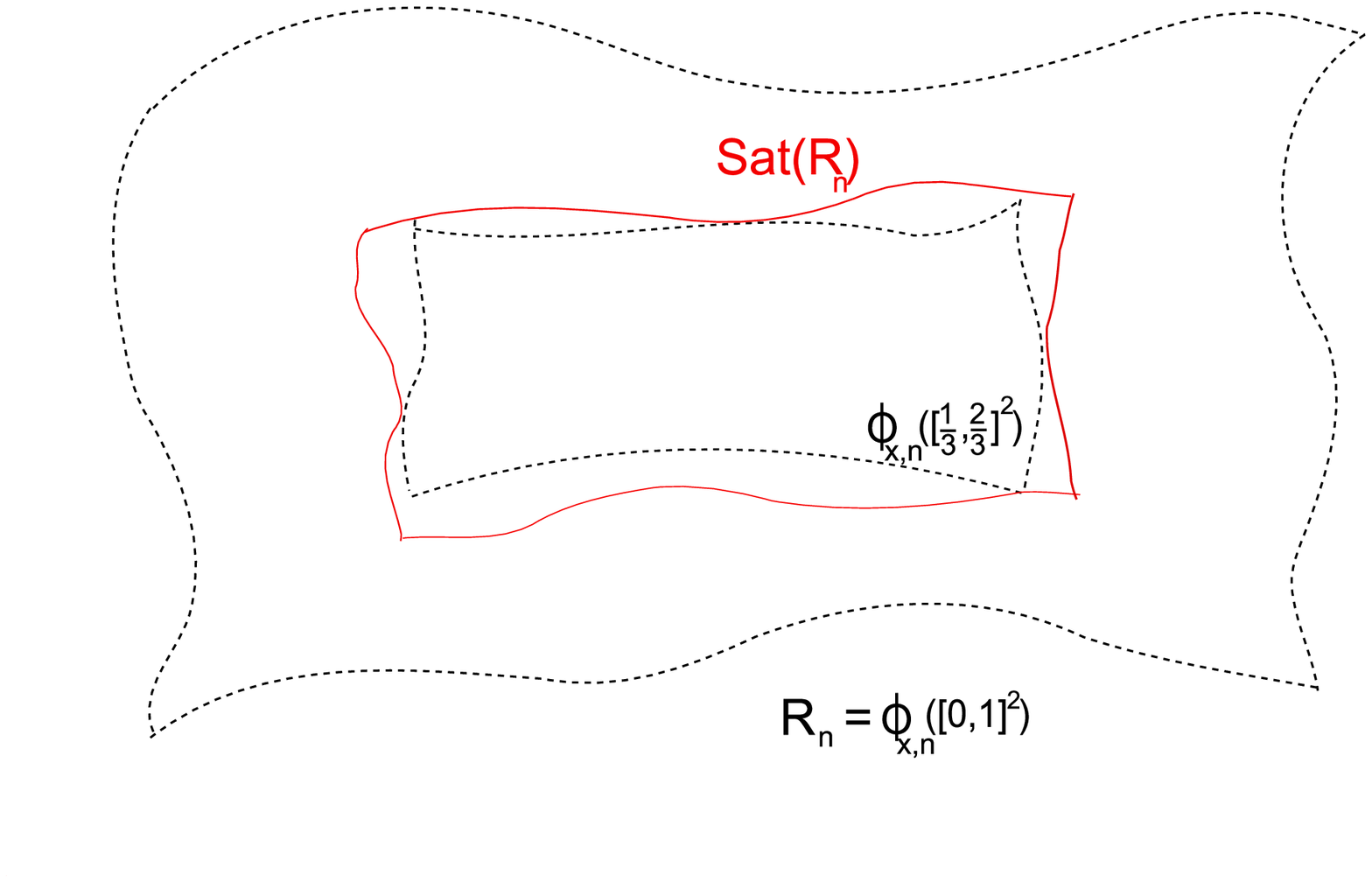}}
\end{prop}

\begin{demo}
 Let $(t_0,s_0)\in [\frac{1}{3},\frac{2}{3}]^2$. Denote $\alpha(t,s_0)\in[0,1]$ the $\mathcal{C}^1$ map such that 
$\alpha(t_0,s_0)=s_0$ and  $\phi_{x,n}(t,\alpha(t,s_0))= \mathcal{F}_n(\phi_{x,n}(t,\alpha(t,s_0)))\cap \mathcal{E}_{n+1}(\phi_{x,n}(t_0,s_0))$ for 
all $t\in [0,1]$. Observe that $\frac{d}{dt}\phi_{x,n}(t,\alpha(t,s_0))$ is colinear to $e_{n+1}(\phi_{x,n}(t,\alpha(t,s_0)))$.

\scalebox{0.38}{\includegraphics{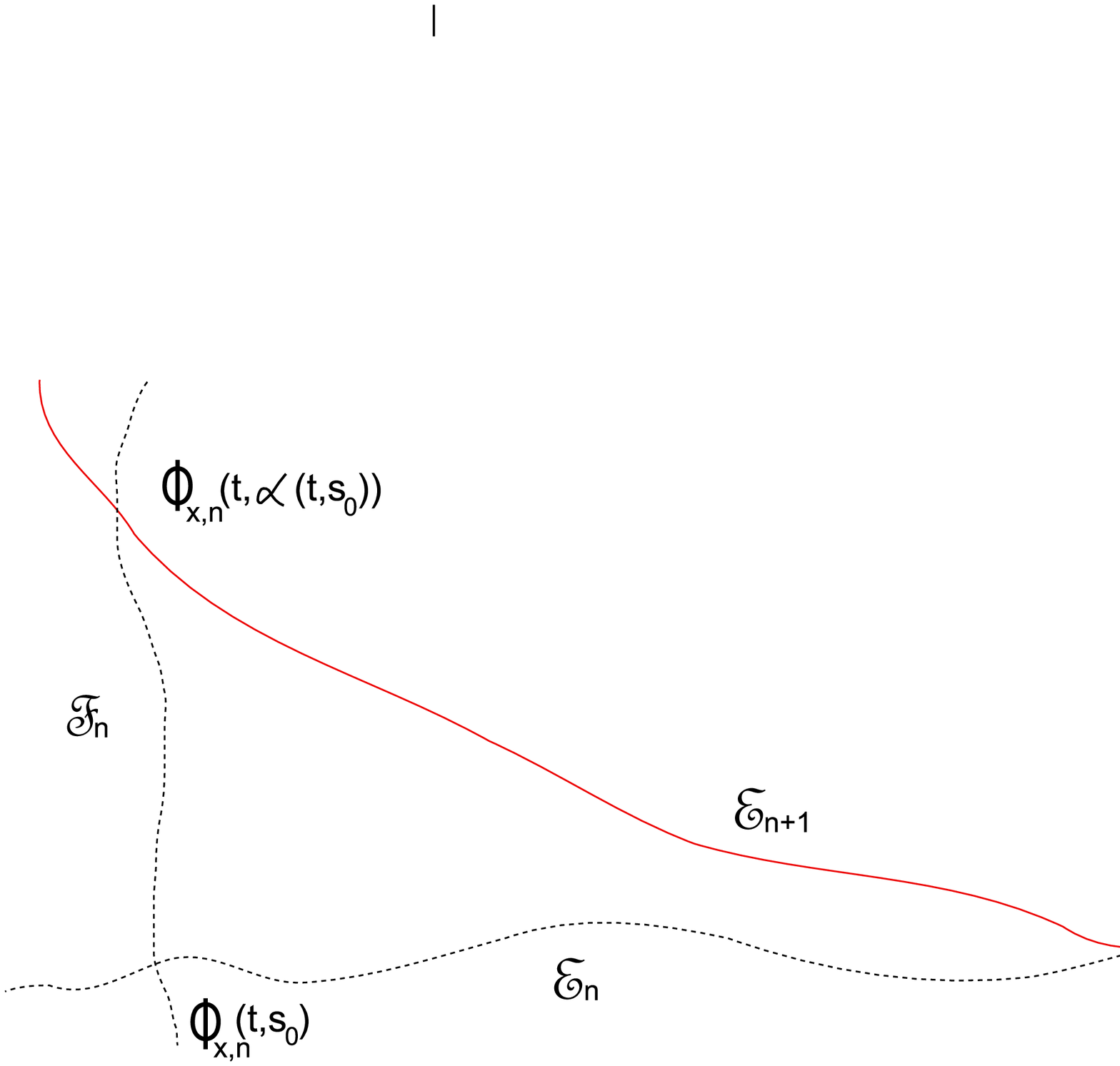}}

But 

\begin{eqnarray*}
\frac{d}{dt}\phi_{x,n}(t,\alpha(t,s_0))&=&\partial_1\phi_{x,n}(t,\alpha(t,s_0))+\partial_1\alpha(t,s_0)\partial_2\phi_{x,n}(t,\alpha(t,s_0))
\end{eqnarray*} 

and then according to Lemma \ref{holangle} and Proposition \ref{angle} we have with $z:=\phi_{x,n}(t,\alpha(t,s_0))$  :
$$(1-K)\frac{l(\mathcal{F}_n(x)\cap R_n)}{l(\mathcal{E}_n(x)\cap R_n)}|\partial_1\alpha(t,s_0)|\leq \tan|\angle| e_n(z),e_{n+1}(z)\leq \frac{(2+K)S(\mathcal{T})^2}{\|D_xT^n\|\|D_{T^nx}T^{-n}\|}$$ 

It follows from the hyperbolicity assumption and the hypothesis $(G_n)$ that  for $n$ large enough we have $|\partial_1\alpha(t,s_0)|\leq K$  and thus   
$\alpha(t,s_0)\in [\frac{1}{3}-K,\frac{2}{3}+K]$ for all $t\in [0,1]$. It means that $\mathcal{E}_{n+1}\left(\phi_{x,n}(t_0,s_0)\right)\cap R_n\subset \phi_{x,n}\left([0,1]\times [\frac{1}{3}-K,\frac{2}{3}+K]\right)$. Similarly we prove that $\mathcal{F}_{n+1}\left(\phi_{x,n}(t_0,s_0)\right)\cap R_n\subset \phi_{x,n}\left([\frac{1}{3}-K,\frac{2}{3}+K]\times [0,1]\right)$.

\end{demo}

We show now that the conclusions of Proposition \ref{marre}  and \ref{angle} apply to the saturated $n+1$-rectangle of a $n$-rectangle satisfying the assumptions $(H_n)$, $(H_{n+1})$ and $(G_n)$. 

\begin{prop}\label{encore}
For any $n$-rectangle $(\phi_{x,n},R_n)$ with $n$ large included in 
$\mathcal{H}_{\mathcal{T}}^n(\chi^+,\chi^-,\gamma,C)$ satisfying $(H_n)$, $(H_{n+1})$ and $(G_n)$ the saturated $n+1$-rectangle  $Sat(R_n)$ satisfies  for all $x'\in R_n$ and $y'\in Sat(R_n)$ :

\begin{eqnarray}\label{nouveau}
l(\mathcal{F}_{n+1}(y')\cap Sat(R_n))& \simeq^{1+K} & \frac{l(\mathcal{F}_n(x')\cap R_n)}{3}\\
l(\mathcal{E}_{n+1}(y')\cap Sat(R_n))& \simeq^{1+K} & \frac{l(\mathcal{E}_n(x')\cap R_n)}{3}
\end{eqnarray}
\end{prop}

\begin{demo}
The lemma follows easily from the following relations which hold for large $n$ :
\begin{itemize}
\item $\forall z,z'\in R_n, \ \angle e_{n+1}(z),e_{n+1}(z')< K $ by Proposition \ref{angles} ;
\item $\forall z \in R_n, \ \angle e_n(z),e_{n+1}(z)<K$ by Lemma \ref{holangle} ;
\item $\phi_{x,n}([\frac{1}{3},\frac{2}{3}]^2)\subset Sat(R_{n})\subset \phi_{x,n}([\frac{1}{3}-K,\frac{2}{3}+K]^2)$ by Proposition \ref{sat}.
\end{itemize}
\end{demo}

\begin{prop}\label{best}
For any $n$-rectangle $(\phi_{x,n},R_n)$ with $n$ large included in 
$\mathcal{H}_{\mathcal{T}}^n(\chi^+,\chi^-,\gamma,C)$ satisfying $(H_n)$, $(H_{n+1})$ and $(G_n)$ and for any $y\in Sat(R_n)$, the saturated $n+1$-rectangle  $(\phi_{y,n+1},Sat(R_n))$  has the following properties for all $t,s\in [0,1]^2$ :

 $$\|\partial_1\phi_{y,n+1}(t,s)\|  \simeq^{1+K} l(\mathcal{E}_{n+1}(y)\cap Sat(R_n))$$
 $$\|\partial_2\phi_{y,n+1}(t,s)\|   \simeq^{1+K}  l(\mathcal{F}_{n+1}(y)\cap Sat(R_n))$$
\end{prop}

The proof of Proposition \ref{best} is technical and may be skipped in a first reading.

\begin{demo}
Let us first prove that $\|\partial_1\phi_{y,n+1}(t,s)\|  \simeq^{1+K} l(\mathcal{E}_{n+1}(y)\cap R_n)$ for all $(t,s)\in [0,1]^2$. It follows from Proposition \ref{angles}  that $\|D(f_{n+1}\circ\phi_{x,n})\|\leq K$. Then by Proposition \ref{angle}, we have for all $z\in R_n$ :
$$\|D_zf_{n+1}\|\leq K l(\mathcal{F}_n(x)\cap R_n)^{-1}$$ 
Let $z$ and $z'$ be two points of  $Sat(R_n)$ lying on the same $n+1$ stable manifold. We reparametrize $\mathcal{F}_{n+1}(z)\cap Sat(R_n)$ and $\mathcal{F}_{n+1}(z')\cap Sat(R_n)$ by arclength. We get two curves 
$z_{n+1}:[-S(z),T(z)]\rightarrow \mathcal{F}_{n+1}(z)$ and $z'_{n+1}:[-S(z'),T(z')]\rightarrow \mathcal{F}_{n+1}(z')$ with $z_{n+1}(0)=z$ and $z'_{n+1}(0)=z'$. We have for all $-\min(S(z),S(z'))\leq u \leq \min(T(z),T(z'))$ :
\begin{eqnarray*}
\|z_{n+1}(u)-z'_{n+1}(u)\|&\leq &\|z-z'\| +\int_{0}^t\|f_{n+1}(z_{n+1}(s))-f_{n+1}(z'_{n+1}(s))\|ds\\
&\leq & \|z-z'\| +\int_{0}^t\|D_{\phi_{x,n}}f_{n+1}\|\|z_{n+1}(s)-z'_{n+1}(s)\|ds
\end{eqnarray*}

By Grönwall Lemma we have then  for all  $-\min(S(z),S(z'))\leq u \leq \min(T(z),T(z'))$ (in particular $|u|\leq l(\mathcal{F}_{n+1}(z)\cap Sat(R_n))\leq l(\mathcal{F}_n(x)\cap R_n)$) :
\begin{eqnarray*}
\|z_{n+1}(u)-z'_{n+1}(u)\|&\leq &\|z-z'\|e^{|u|\|D_{\phi_{x,n}}f_{n+1}\|}\\
&\leq & \|z-z'\|e^{|u|Kl(\mathcal{F}_n(x)\cap R_n)^{-1}}\\
&\leq & (1+K)\|z-z'\|
\end{eqnarray*}

Let us denote by $Z'_{n+1}(u)\in Sat(R_n)$ the intersection point of $\mathcal{E}_{n+1}(z_{n+1}(u))$ and $\mathcal{F}_{n+1}(z')$. Remark that by continuity of the $n+1$ stable field $e_{n+1}$ there exists $w$ in the finite time unstable piece $[z_{n+1}(u),Z'_{n+1}(u)]$  such that the vector $z_{n+1}(u)-Z'_{n+1}(u)$ is colinear to $e_{n+1}(w)$. Then since $|\angle| e_{n+1}(v),e_{n+1}(v')\leq K$ for all $v,v'\in R_n$ it follows that $\mathcal{F}_{n+1}(z)\cap Sat(R_n)$ (resp. $\mathcal{F}_{n+1}(z')\cap Sat(R_n)$)  lies in the cone $\{v, |\angle| v-z_{n+1}(u),f_{n+1}(w)<K\}$ (resp.   $\{v, |\angle| v-Z'_{n+1}(u),f_{n+1}(w)<K\}$). In particular the distance $d(z_{n+1}(u),\mathcal{F}_{n+1}(z'))$ is bounded from below by $|z_{n+1}(u)-Z'_{n+1}(u)|\cos K$. It follows that :

\begin{eqnarray*}
\|z_{n+1}(u)-Z'_{n+1}(u)\|&\leq &  \frac{d(z_{n+1}(u),\mathcal{F}_{n+1}(z))}{\cos K}\\
&\leq & (1+K)\|z_{n+1}(u)-z'_{n+1}(u)\|\\
&\leq &(1+K)\|z-z'\|
\end{eqnarray*}

\scalebox{0.45}{\includegraphics{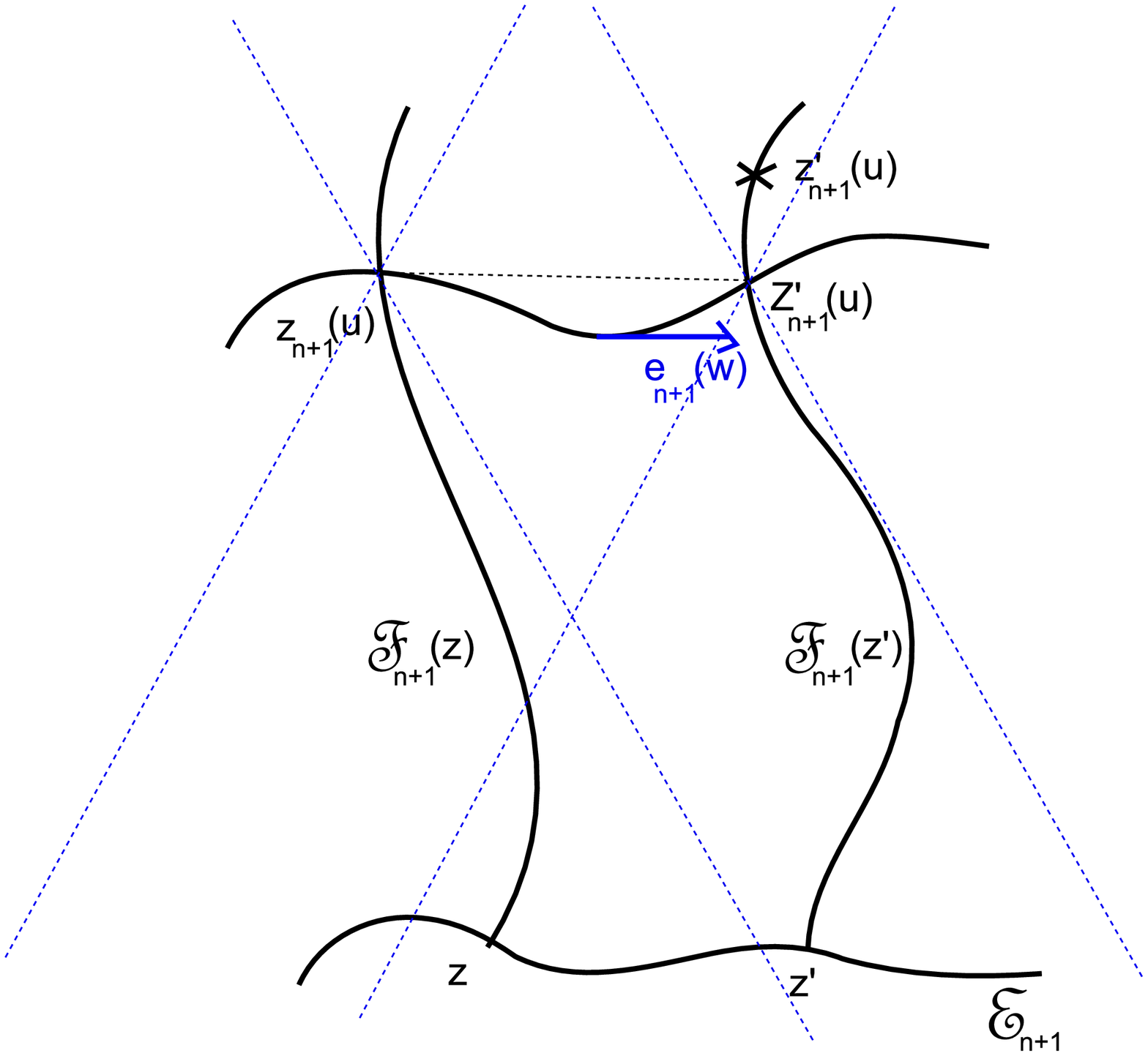}}

Let $t,t',s,s'$ such that $\phi_{y,n+1}(t,s)=z$, $\phi_{y,n+1}(t',s)=z'$, $\phi_{y,n+1}(t,s')=z_{n+1}(u)$ and $\phi_{y,n+1}(t',s')=Z'_{n+1}(u)$. When $z$ goes to $z'$ we get from the last inequality :
$\|\partial_1 \phi_{y,n+1}(t,s')\|\leq (1+K)\|\partial_1 \phi_{y,n+1}(t,s)\|$ and by inverting the roles of $s$ and $s'$, we have for all $t,s,s' \in [0,1]$ :
$$\|\partial_1 \phi_{y,n+1}(t,s')\|\simeq^{1+K}\|\partial_1 \phi_{y,n+1}(t,s)\|$$

But the $n+1$ stable manifold of $y$ is reparametrized with constant rate. It follows that :

 $$\|\partial_1 \phi_{y,n+1}(t,s')\|\simeq^{1+K}l(\mathcal{E}_{n+1}(y)\cap Sat(R_n))$$

We deal now with the partial derivative $\partial_2 \phi_{y,n+1}$. By Proposition \ref{angles} applied  to $T^{n+1}\circ \phi_{x,n}$ for $T^{-1}$, we have :
\begin{eqnarray}\label{mdr}
\|D(e_{n+1}^{(n+1)}\circ \phi_{x,n})\|\leq K
\end{eqnarray}

Let us compute $D(T^{n+1}\circ \phi_{x,n})^{-1}=D(\phi_{x,n}^{-1})DT^{-n-1}$. 
In the bases $(e_{n+1}^{(n+1)},f_{n+1}^{(n+1)})\rightarrow (e_{n+1},f_{n+1})$ the matrix associated to $DT^{-n-1}$ is given by : $$\begin{pmatrix}
\|(DT^{n+1})^{-1}\| & 0\\
0& \|DT^{n+1}\|^{-1} 
\end{pmatrix}$$

and by Proposition \ref{angle} the matrix associated to $D(\phi_{x,n}^{-1})$ in the bases $(e_{n},f_{n})\rightarrow (\partial_t,\partial_s)$ is up to a factor $1+K$ : 
$$\begin{pmatrix}
l(\mathcal{E}_n(x)\cap R_n)^{-1} & 0\\
0&l(\mathcal{F}_n(x)\cap R_n)^{-1} 
\end{pmatrix}$$

Finally the matrix passage from the basis $(e_{n+1},f_{n+1})$ to the basis $(e_{n},f_{n})$ is just the rotation matrix with angle $\angle e_{n+1},e_n$. After a simple computation we get for all $z\in T^{n+1}R_n$ up to a factor $1+K$ :

$$D_z(T^{n+1}\circ \phi_{x,n})^{-1} = 
 \begin{pmatrix}
\frac{\cos(\angle e_{n+1},e_n)\|(D_xT^{n+1})^{-1}\|}{l(\mathcal{E}_n(x)\cap R_n)} & -\frac{\sin(\angle e_{n+1},e_n)}{\|D_xT^{n+1}\|l(\mathcal{E}_n(x)\cap R_n)}\\
\frac{\sin(\angle e_{n+1},e_n)\|(D_xT^{n+1})^{-1}\|}{l(\mathcal{F}_n(x)\cap R_n)}& \frac{\cos(\angle e_{n+1},e_n)}{\|D_xT^{n+1}\|l(\mathcal{F}_n(x)\cap R_n)} 
\end{pmatrix}$$

 We deduce from the inequalities $\tan|\angle| e_n(x),e_{n+1}(x)<\frac{2S(\mathcal{T})^2}{\|(D_xT^n)^{-1}\|\|D_xT^n\|}$ (Lemma \ref{holangle}) and 
$(1+K)l(\mathcal{E}_n(x)\cap R_n)\geq l(\mathcal{F}_n(x)\cap R_n) \geq (1-K)\frac{l(\mathcal{E}_n(x)\cap R_n)}{\|D_xT^n\|}$ (Hypothesis $(G_n)$) that for $n$ large enough :

$$\|D_z(T^{n+1}\circ \phi_{x,n})^{-1}\| \leq (1+K)l(\mathcal{E}_n(x)\cap R_n)^{-1}\|(D_xT^{n+1})^{-1}\|$$

and by Proposition \ref{encore} and assumption $(H_{n+1})$ :

\begin{eqnarray*} \|D_z(T^{n+1}\circ \phi_{x,n})^{-1}\|&\leq & \frac{1+K}{3}l(\mathcal{E}_{n+1}(y)\cap Sat(R_{n}))^{-1}\|(D_{x}T^{n+1})^{-1}\|\\
&\leq & l\left(T^{n+1}\mathcal{E}_{n+1}(y)\cap T^{n+1}Sat(R_{n})\right)^{-1}
\end{eqnarray*}

Then we have by Equation (\ref{mdr}) : $$\|D_z\left(e^{(n+1)}_{n+1}\circ T^{-n-1}\right)\|\leq Kl\left(T^{n+1}\mathcal{E}_{n+1}(y)\cap T^{n+1}Sat(R_{n})\right)^{-1}$$ and thus we can apply the first part of the proof to the $n+1$-rectangle  $T^{n+1}Sat(R_{n})$ for $T^{-1}$ to get :

$$\|\partial_2(T^{n+1}\circ \phi_{y,n+1})\|\simeq^{1+K} l\left(T^{n+1}\mathcal{F}_{n+1}(y)\cap T^{n+1}Sat(R_n)\right)$$

and therefore again by assumption $(H_{n+1})$ :

$$\|\partial_2\phi_{y,n+1}\| \simeq^{1+K}  l(\mathcal{F}_{n+1}(y)\cap Sat(R_n))$$
\end{demo}

We deduce from the previous proposition that the change of coordinates $\phi_{x,n}^{-1}\circ\phi_{y,n+1}$ from the $n+1$ to the $n$ foliations is $\mathcal{C}^1$ bounded. 

\begin{coro}\label{trois}
For any $n$-rectangle $(\phi_{x,n},R_n)$ with $n$ large included in 
$\mathcal{H}_{\mathcal{T}}^n(\chi^+,\chi^-,\gamma,C)$ satisfying $(H_n)$, $(H_{n+1})$ and $(G_n)$ and for any $y\in Sat(R_n)$, the saturated $n+1$-rectangle  $(\phi_{y,n+1},Sat(R_n))$ satisfies :\\

\begin{center}
\begin{itemize}
\item $\|D\left((\phi_{x,n})^{-1}\circ\phi_{y,n+1}\right)\|\leq \frac{1}{3}+K$ ;\\
\item $\|D\left((\phi_{y,n+1})^{-1}\circ(\phi_{x,n})_{/[\frac{1}{3},\frac{2}{3}]^2}\right)\|\leq 3+K$.
\end{itemize}
\end{center}
\end{coro}

\begin{demo}

The previous propostion claims that 
in the bases $(\partial_t,\partial_s)\rightarrow (e_{n+1},f_{n+1})$ the matrix associated to $D\phi_{y,n+1}$ is  up to a factor $1+K$ : 

$$\begin{pmatrix}
l(\mathcal{E}_{n+1}(y)\cap Sat(R_n)) & 0\\
0&l(\mathcal{F}_{n+1}(y)\cap Sat(R_n)) 
\end{pmatrix}$$

then as in the previous proof we get after an easy matricial computation that for all $z\in [0,1]^2$ and up to a factor $1+K$ :

$$D_z\left((\phi_{x,n})^{-1}\circ \phi_{y,n+1}\right) = 
 \begin{pmatrix}
\cos(\angle e_{n+1},e_n)\frac{l(\mathcal{E}_{n+1}(x)\cap Sat(R_n))}{l(\mathcal{E}_n(x)\cap R_n)} & \sin(\angle e_{n+1},e_n)\frac{l(\mathcal{F}_{n+1}(x)\cap Sat(R_n))}{l(\mathcal{E}_n(x)\cap R_n)}\\
-\sin(\angle e_{n+1},e_n)\frac{l(\mathcal{E}_{n+1}(x)\cap Sat(R_n))}{l(\mathcal{F}_n(x)\cap R_n)}&\cos(\angle e_{n+1},e_n)\frac{l(\mathcal{F}_{n+1}(x)\cap Sat(R_n))}{l(\mathcal{F}_n(x)\cap R_n)} 
\end{pmatrix}$$

Then we conclude the proof by Proposition \ref{encore}, assumption $(G_n)$ and Lemma \ref{holangle}.
\end{demo}

\section{Reparametrization of $(n,\epsilon)$ Bowen balls}
This section is devoted to the proof of Proposition \ref{main}. We begin with two technical lemmas.

\subsection{Combinatorial Lemma}

We recall first a usual combinatorial lemma which was already used by T.Downarowicz and A.Maass in \cite{Dow}.

\begin{defi}
Let $S\in \N$ and $n\in\N$. We say that a sequence of $n$ positive integers $\mathcal{K}_n:=(k_1,...,k_n)$ admits the value $S$ if $\frac{1}{n}\sum_{i=1}^nk_i\leq S$. 
\end{defi}

The number of  sequences of $n$ positive integers admitting the value $S$ is exactly the binomial coefficient ${nS \choose n}$. By applying Stirling's formula we get easily :   

\begin{lem}\label{combi}
The logarithm of the number of sequences of $n$ positive integers admitting the value $S$ is at most 
$nSH(S)+1$.
\end{lem}

\subsection{Lemma of Linear algebra}

The following basic algebraic lemma allows us to control the oscillation of $D_{T^n\circ\phi}T^{-n}$ and $D_{\phi}T^n$ at the same time. It claims that the defect of multiplicativity of the norm of a product of two matrices of $\mathcal{GL}_2(\R)$ is almost equal to this  of the product of the inverses.

\begin{lem}\label{alg}
Let $A,B\in\mathcal{GL}_2(\R)$ then  
$$\frac{\|AB\|}{\|A\| \|B\|} =\frac{\|B^{-1}A^{-1}\|}{\|A^{-1}\|\|B^{-1}\|}$$
\end{lem}

\begin{demo}
This follows immediately from the formula $det(C) =\frac{\|C\|}{\|C^{-1}\|}$ which holds for all matrices $C$ of $\mathcal{GL}_2(\R)$ and from the multiplicativity of the determinant. 
\end{demo}

This fact is specific to the dimensions $1$ and $2$. Indeed consider for $x\neq 0$ the matrices $A(x),B(x)\in \mathcal{GL}_3(\R)$ defined by :

\begin{eqnarray*}
A(x)=\begin{pmatrix}
1 & 0 & 0 \\
0 & x & -1 \\
0 & x & 1
\end{pmatrix} & and & B(x)=\begin{pmatrix}
x & x & 0 \\
1 & -1 & 0 \\
0 & 0 & 1   
\end{pmatrix}
\end{eqnarray*}

One easily computes : 

\begin{eqnarray*}
A(x)^{-1}=\begin{pmatrix}
1 & 0 & 0 \\
0 & \frac{1}{2x} & \frac{1}{2x} \\
0 & -\frac{1}{2} & \frac{1}{2}
\end{pmatrix} & and & B(x)^{-1}=\begin{pmatrix}
\frac{1}{2x} & \frac{1}{2} & 0 \\
\frac{1}{2x} & \frac{-1}{2} & 0 \\
0 & 0 & 1
\end{pmatrix}
\end{eqnarray*}

\begin{eqnarray*}
A(x)B(x)=\begin{pmatrix}
x & x & 0 \\
x & -x & -1 \\
x & -x & 1
\end{pmatrix}& and & B(x)^{-1}A(x)^{-1}=\begin{pmatrix}
\frac{1}{2x} & \frac{1}{4x} &  \frac{1}{4x}\\
\frac{1}{2x} & \frac{-1}{4x} &  \frac{-1}{4x}\\
0 & -\frac{1}{2} & \frac{1}{2}
\end{pmatrix}
\end{eqnarray*}

In particular we have $\|A(x)\|_F=\|B(x)\|_F \sim \sqrt{2}x$, $\|A(x)^{-1}\|_F=\|B(x)^{-1}\|_F \sim \sqrt{\frac{3}{2}}$, $\|A(x)B(x)\|_F\sim \sqrt{6}x$ but $\|B(x)^{-1}A(x)^{-1}\|_F\sim \sqrt{\frac{1}{2}}$ when $x$ goes to $+\infty$.\\

 We give now a dynamical interpretation of the previous lemma in our setting. In the proof of Theorem \ref{gg} assuming Proposition 1 the quantity  $\lambda^+_n(x,T)-\frac{1}{n}\log \|D_xT^{n}\|$ is close to $\overline{\chi^+}(\mu)-\chi^{+}(\nu)$ when $n$ goes to infinity for typical $\nu$ points $x$ (where $\nu$ is an ergodic measure near an invariant measure $\mu$). If we define similarly   $\lambda^-_n(x,T)=\frac{1}{n}\sum_{l=0}^{n-1}\log^+\|(D_{T^lx}T)^{-1}\|$ then we can also ensure that 
$\lambda^-_n(x,T)-\frac{1}{n}\log \|(D_xT^{n})^{-1}\|$ is close to $-\overline{\chi^-}(\mu)+\chi^{-}(\nu)$ when $n$ goes to infinity. Then by applying the previous algebraic lemma at each time and by averaging along the orbit of $x$ we obtain  that  $\overline{\chi^+}(\mu)-\chi^{+}(\nu)\simeq -\overline{\chi^-}(\mu)+\chi^{-}(\nu)$ for $\nu$ close to $\mu$. But observe now this fact follows from the continuity in $\mu\in \mathcal{M}(X,T)$ of $\overline{\chi^{+}+\chi^{-}}(\mu)=\int_{M}\log Jac_xTd\mu(x)$.

\subsection{Proof of Proposition \ref{main}}

As in Yomdin's theory, we consider the local dynamic at one point. We fix a riemanian structure $\|\|$ on the surface $M$ and we denote $d$ the induced distance on $M$, $R_{inj}$ the radius of injectivity and $exp:TM(R_{inj})\rightarrow M$ the exponential map, where $TM(r):=\{(x,u), \ u\in T_xM, \ \|u\|_x<r\}$. There exist
$R<R'<R_{inj}$ such that $T(B(x,R))\subset B(Tx,R')$ for all $x\in M$.\\

 Let $x\in M$ and $n\in \N$. We consider the map 
$T_n^x:T_{T^{n-1}x}M(R)\rightarrow T_{T^nx}M(R')$  defined by $T_n^x=exp^{-1}_{T^nx}\circ T\circ exp_{T^{n-1}x}$. For all $\epsilon<\frac{R}{2}$, we put $T_{n,\epsilon}^x=\epsilon^{-1}T_n^x(\epsilon .):B(0,2)\rightarrow T_{T^nx}M\simeq \R^2$ and  $\mathcal{T}_{\epsilon}^x:=(T_{n,\epsilon}^x)_{n\in \N}$.\\

We choose $\epsilon>0$ be small enough such that $\|D^2T_{n,\epsilon}^x\|\leq \inf_{z\in B(0,2)}\|D_zT_{n,\epsilon}^x\|$ and  $\|D^2(T_{n,\epsilon}^x)^{-1}\|\leq \frac{\inf_{z\in B(0,2)}\|(D_zT_{n,\epsilon}^x)^{-1}\|}{\|DT_{n,\epsilon}^x\|}$ for all integers $n$ and all $x\in M$. One can also assume  that $\|D_xT\|\simeq^{1+K}\|D_yT\|$ and $\|D_{Tx}T^{-1}\|\simeq^{1+K}\|D_{Ty}T^{-1}\|$ for all $x,y\in M$ with  $d(x,y)<2\epsilon$ and therefore $\|D_zT_{n,\epsilon}^x\|\simeq^{1+K}\|D_{z'}T_{n,\epsilon}^x\|$ and  $\|D_{T_{n,\epsilon}^xz}(T_{n,\epsilon}^x)^{-1}\|\simeq^{1+K}\|D_{T_{n,\epsilon}^xz'}(T_{n,\epsilon}^x)^{-1}\|$   for all $n\in\N$ and  for all $z,z'\in B(0,2)$. \\

We define now subsets of the Bowen ball $B(n+1,1)$ where the defect of multiplicativity of the norm of the composition $DT_{i+1}\circ DT^i$ is prescribed at each step $1\leq i\leq n$. 

\begin{defi} Let $\mathcal{T}:=(T_n)_{n\in\N}$ be a sequence of $\mathcal{C}^2$ maps from $B(0,2)$ to $\R^2$ and  let $\mathcal{K}_n:=(k_1,...,k_n)$  be a sequence of $n$ positive integers, we denote $\mathcal{H}(\mathcal{K}_n)$ the subset of $B(n+1,1)$ defined by :

$$\mathcal{H}(\mathcal{K}_n):=\left\{y\in B(n+1,1)\ : \ \forall 1\leq i\leq n, \ \left[ \log^+ \frac{\|D_yT^i\|\max(\|D_{T^iy}T_{i+1}\|,1)}{\|D_yT^{i+1}\|}\right] +1=k_i \right\}$$

\end{defi}

 In the next lemma we estimate the numbers of such sets intersecting a finite time hyperbolic set by using the combinatorial argument of Lemma \ref{combi}.  We write  $\lambda^+_n:=\frac{1}{n}\sum_{i=0}^{n-1}\log^+\|D_{0}T_i\|$. When $\mathcal{T}=\mathcal{T}_{\epsilon}^x$  we have $\lambda_n^+=\lambda_n^+(x,T)=\frac{1}{n}\sum_{i=0}^{n-1}\log^+\|D_{T^ix}T\|$. 

\begin{lem}\label{com}
Let $\chi^+>0>\chi^-$, $\frac{\min(\chi^+,-\chi^-,1)}{3}>\gamma>0$ and $C>1$  and let $\mathcal{T}:=(T_n)_{n\in\N}$ be a sequence of $\mathcal{C}^2$ maps from $B(0,2)$ to $\R^2$ such that $\|D_zT_n\|\simeq^{2}\|D_{z'}T_n\|$ for all $n\in\N$ and  $z,z'\in B(0,1)$. There exists an integer $N$ depending only on $C$ such that for $n>N$ the number of sequences $\mathcal{K}_{n-1}$ such that $\mathcal{H}(\mathcal{K}_{n-1})$ has a non empty intersection with $\mathcal{H}^n_{\mathcal{T}}(\chi_+,\chi_-,\gamma,C)\cap B(n,1)$ is bounded above  by 

$$e^{3n-2}e^{(n-1)(\lambda^+_n-\chi^+)H([\lambda^+_n-\chi^+]+3)}$$

\end{lem}

Indeed if $y\in \mathcal{H}^n_{\mathcal{T}}(\chi_+,\chi_-,\gamma,C)\cap B(n,1)$, then $\sum_{i=0}^{n-1}\log^+\|D_{T^iy}T_{i+1}\|-\log \|D_yT^n\|\leq 
n(\lambda^+_n+\log 2-\chi^++\gamma)+\log C$. Thus the sequence $\left(\left[ \log \frac{\|D_yT^i\|\max(\|D_{T^iy}T_{i+1}\|,1)}{\|D_yT^{i+1}\|}\right]+1\right)_{i=1,...,n-1}$ admits the value  $[\lambda^+_n-\chi^+]+3$ for $n>\frac{\log C}{1-\log 2-\frac{1}{3}}$.  We apply finally the combinatorial Lemma \ref{combi} to conclude the proof of Lemma \ref{com}. \\

Proposition 1 follows now direcly from Lemma \ref{com} and the following Proposition \ref{rep} applied to the sequences $\mathcal{T}=\mathcal{T}_{\epsilon}^x$ for $\epsilon>0$ chosen as previously and for all $x\in M$.

\begin{prop}\label{rep}
Let    $\chi^+>0>\chi^-$, $\frac{\min(\chi^+,-\chi^-,1)}{3}>\gamma>0$ and $C>1$  and  let $\mathcal{T}:=(T_n)_{n\in\N}$ be a 
sequence of $\mathcal{C}^2$ maps from $B(0,2)$ to $\R^2$ with $S(\mathcal{T})<+\infty$ and such that  $\|D^2T_n\|\leq \inf_{z\in B(0,2)}\|D_zT_n\|$, 
$\|D^2T_n^{-1}\|\leq \frac{\inf_{z\in B(0,2)}\|(D_zT_n)^{-1}\|}{\|DT_n\|}$, $\|D_zT_n\|\simeq^{1+K}\|D_{z'}T_n\|$, $\|D_{T_nz}T_n^{-1}\|\simeq^{1+K}\|D_{T_nz'}T_n^{-1}\|$ for all 
integers $n$ and for all $z,z'\in B(0,2)$. Then there exists a real number $D$ depending only on $\chi^+,\chi^-,\gamma,C$ and $S(\mathcal{T})$ and an universal constant $A$ such 
that for all sequences $\mathcal{K}_{n-1}=(k_1,...,k_{n-1})$ of  
$n-1$ positive integers there exists a family $\mathcal{F}_n$ of admissible charts of  $n$-rectangles included in  $\mathcal{H}^n_{\mathcal{T}}(\chi^+,\chi^-,\gamma,2C)$ 
satisfying :

\begin{enumerate}[(i)]
\item  $\forall \phi_n\in \mathcal{F}_n, \ \phi_n([0,1]^2)\subset B(n+1,2)$ ;
\item  $\forall \phi_n\in \mathcal{F}_n \ \|D\phi_n\|\leq \frac{1}{4}$ and $\forall 1\leq k\leq n, \ \|D(T^{k}\circ \phi_n)\|\leq 1$ ;  
\item  $\forall \phi_n\in \mathcal{F}_n \ \forall 0\leq k\leq n, \ \|D(x\mapsto D_{\phi_n(x)}T^{k})\|\leq K\|D_{\phi_n}T^k\|$ ;
\item  $\forall \phi_n\in \mathcal{F}_n \ \forall 0 \leq k\leq n, \ \|D(x\mapsto D_{T^k\circ\phi_n(x)}T^{-k})\|\leq K\|D_{T^k\circ\phi_n}T^{-k}\|$ ;
\item  $\forall \phi_n\in \mathcal{F}_n \ \forall x\in [0,1]^2,  \ (1+K)l(\mathcal{E}_{n}(x)\cap \phi_n([0,1]^2))\geq l(\mathcal{F}_{n}(x)\cap \phi_n([0,1]^2)) \geq (1-K)\frac{l(\mathcal{E}_{n}(x)\cap \phi_n([0,1]^2))}{\max_{k=0,...,n}\|D_xT^k\|}$ ;
\item  $\mathcal{H}^n_{\mathcal{T}}(\chi^+,\chi^-,\gamma,C)\cap \mathcal{H}(\mathcal{K}_{n-1})\cap B(n+1,1)\subset \bigcup_{\phi_n\in \mathcal{F}_n}\phi_n([\frac{4}{9},\frac{5}{9}]^2)$ ;
\item  $\log \sharp \mathcal{F}_n \leq D+An+2\sum_{i=1}^{n-1}k_i$.
\end{enumerate}
\end{prop}

With the previous notations the assumptions (iii)+(iv) and  (v) are respectively  $(H_{k})$ for all $0\leq k\leq n$ and $(G_{n})$. 

\begin{demo}
We argue by induction on $n\geq N$ where $N$ is an integer such that all the statements of Subsection \ref{supsub} holds for $n>N$. The initial step is easily checked. We assume now the existence of the family $\mathcal{F}_n$ and 
we built $\mathcal{F}_{n+1}$. 
Let $\phi_n\in \mathcal{F}_n$.\\

\underline{\textbf{First step :} We divide the $n$-rectangle $\phi_n$ in $n$-subrectangles $\phi'_n$ satisfying $(H_{n+1})$.}\\
 
We cover the square $[\frac{4}{9},\frac{5}{9}]^2$ into $([K^{-1}e^{k_n}]+1)^2$ subsquares  of size $\frac{1}{[K^{-1}e^{k_n}]+1}<Ke^{-k_n}$ such that by reparametrizing these subsquares from the unit square by affine contractions $\psi$  we have $\bigcup_{\psi}\psi([\frac{4}{9},\frac{5}{9}]^2)=[\frac{4}{9},\frac{5}{9}]^2$. By composing $\psi$ with $\phi_{n}$ we get new  $n$-rectangles $\phi'_{n}=\phi_n\circ \psi$ such that   
$\bigcup_{\phi'_n}\phi'_{n}([\frac{4}{9},\frac{5}{9}]^2)=\bigcup_{\phi_n}\phi_n([\frac{4}{9},\frac{5}{9}]^2)$. Note that by Proposition \ref{angles}   the $n$-rectangles $\phi'_{n}$ satisfy again the assumption $(G_n)$. From now on we only  consider $n$-rectangles $\phi'_{n}$ whose middle ninth intersect the set $\mathcal{H}^{n+1}_T(\chi_+,\chi_-,\gamma,C)\cap\mathcal{H}(\mathcal{K}_n)\cap B(n+1,1)$. Fix such a $n$-rectangle $\phi'_n$ and choose $w\in [\frac{4}{9},\frac{5}{9}]^2$ such that $\phi'_n(w)$ belongs to this set. Let us show that $\|D(x\mapsto 
D_{\phi'_n(x)}T^{n+1})\|\leq K\|D_{\phi'_{n}}T^{n+1}\|$. We compute for all $y\in [0,1]^2$ :

\begin{multline*}
D_y(x\mapsto D_{\phi'_{n}(x)}T^{n+1})=D_y(x\mapsto D_{T^n\circ \phi'_{n}(x)}T_{n+1})D_{\phi'_{n}(y)}T^{n} \\
 +D_{T^n\circ\phi'_{n}(y)}T_{n+1}D_y(x\mapsto D_{\phi'_{n}(x)}T^{n})  
\end{multline*}

We bound the norm of the first term $D_y(x\mapsto D_{T^n\circ \phi'_{n}(x)}T_{n+1})$ as follows :
\begin{eqnarray*}
\|D_y(x\mapsto D_{T^n\circ \phi'_{n}(x)}T_{n+1})\|&\leq &\|D^2_{T^{n}\circ \phi'_n}T_{n+1}\| \|D_y(T^n\circ \phi'_{n})\|\\
&\leq & \|D^2_{T^{n}\circ \phi'_n}T_{n+1}\|Ke^{-k_n}\|D_{\psi(y)}(T^n\circ \phi_{n})\|\\
&\leq & Ke^{-k_n}\|D_{T^n\circ\phi'_n(y)}T_{n+1}\|
\end{eqnarray*}
where the last inequality follows from  $\|D^2T_{n+1}\|\leq \inf_{z\in B(0,2)}\|D_zT_{n+1}\|$ and the induction hypothesis (ii). Then we get according to the induction hypothesis (iii) and $\|D_zT_{n+1}\|\simeq^{1+K}\|D_{z'}T_{n+1}\|$ for all $z,z'\in B(0,2)$ :
 
\begin{eqnarray*} 
\|D_y(x\mapsto D_{\phi'_{n}(x)}T^{n+1})\| &\leq &   Ke^{-k_n}\|D_{T^n\circ\phi'_n(y)}T_{n+1}\|\left(\|D_{\phi'_n(y)}T^n\|+\|D_{\psi(y)}(x\mapsto D_{\phi_n(x)}T^n)\|\right)\\
& \leq &  Ke^{-k_n}\|D_{T^n\circ\phi'_n(y)}T_{n+1}\|\left(\|D_{\phi'_n(y)}T^n\|+\|D_{\phi_n}T^n\|\right) \\
 &\leq & Ke^{-k_n}\|D_{T^n\circ\phi'_n(w)}T_{n+1}\|\|D_{\phi'_n(w)}T^n\|
 \end{eqnarray*}
 
and as $\phi'_n(w)\in \mathcal{H}(\mathcal{K}_n)$ we obtain :
  
\begin{eqnarray*}
\|D_y(x\mapsto D_{\phi'_{n}(x)}T^{n+1})\| &\leq & K\|D_{\phi'_n(w)}T^{n+1}\|\\
 &\leq & K\|D_{\phi'_n}T^{n+1}\|
\end{eqnarray*}

Similarly we have :

\begin{multline*}
D_y(x\mapsto D_{T^{n+1}\circ \phi'_n(x)}T^{-(n+1)})=D_{T^n\circ\phi'_{n}(y)}T^{-n}D_y(x\mapsto D_{T^{n+1}\circ \phi'_{n}(x)}T_{n+1}^{-1}) \\
+D_y(x\mapsto D_{T^n\circ \phi'_{n}(x)}T^{-n})D_{T^{n+1}\circ\phi'_{n}(y)}T_{n+1}^{-1}  
\end{multline*}

The norm of the  term $D_y(x\mapsto D_{T^{n+1}\circ \phi'_{n}(x)}T_{n+1}^{-1})$ is bounded from above in the following way 

\begin{eqnarray*}
\|D_y(x\mapsto D_{T^{n+1}\circ \phi'_{n}(x)}T_{n+1}^{-1})\|&\leq &\|D^2_{T^{n+1}\circ \phi'_n}T_{n+1}^{-1}\| \|D_y(T^{n+1}\circ \phi'_{n})\|\\
&\leq & \|D^2_{T^{n+1}\circ \phi'_n}T_{n+1}^{-1}\|\|D_{T^n\circ \phi'_n(y)}T_{n+1}\| \|D_{y}(T^n\circ \phi'_{n})\|\\
&\leq & \|D^2_{T^{n+1}\circ \phi'_n}T_{n+1}^{-1}\|\|D_{T^n\circ \phi'_n(y)}T_{n+1}\| Ke^{-k_n}\|D_{\psi(y)}(T^n\circ \phi_{n})\|\\
&\leq & Ke^{-k_n}\|D_{T^{n+1}\circ\phi'_{n}(y)}T_{n+1}^{-1}\|
\end{eqnarray*}
where the last inequality follows from  $\|D^2T_{n+1}^{-1}\|\leq \frac{\inf_{z\in B(0,2)}\|(D_{z}T_{n+1})^{-1}\|}{\|DT_{n+1}\|}$ and the induction hypothesis (ii). Then we get according to the induction hypothesis (iv) and $\|D_{T_{n+1}z}T_{n+1}^{-1}\|\simeq^{1+K}\|D_{T_{n+1}z'}T_{n+1}^{-1}\|$ for all $z,z'\in B(0,2)$ :

 \begin{eqnarray*}
 \|D(x\mapsto D_{T^{n+1}\circ \phi'_{n}(x)}T^{-(n+1)})\|&\leq & Ke^{-k_n}\|D_{T^{n+1}\circ\phi'_{n}(y)}T_{n+1}^{-1}\|\left(\|D_{T^n\circ\phi'_{n}(y)}T^{-n}\|+\|D_{\psi(y)}(x\mapsto D_{T^n\circ \phi_n(x)}T^{-n})\|\right)\\
 &\leq & Ke^{-k_n} \|D_{T^{n+1}\circ\phi'_{n}(y)}T_{n+1}^{-1}\|\left(\|D_{T^n\circ\phi'_{n}(y)}T^{-n}\|+\|D_{T^n\circ\phi_{n}}T^{-n}\|\right)\\
 &\leq & Ke^{-k_n} \|D_{T^{n+1}\circ\phi'_{n}(w)}T_{n+1}^{-1}\|\|D_{T^n\circ\phi'_{n}(w)}T^{-n}\|
\end{eqnarray*}

Now we have by applying Lemma \ref{alg} : 
$$\frac{\|D_{T^{n+1}\circ\phi'_{n}(w)}T_{n+1}^{-1}\|\|D_{T^n\circ\phi'_{n}(w)}T^{-n}\|}{\|D_{T^{n+1}\circ\phi'_{n}(w)}T^{-(n+1)}\|}= \frac{\|D_{\phi'_n(w)}T^{n}\|\|D_{T^n\circ \phi'_n(w)}T_{n+1}\|}{\|D_{\phi'_n(w)}T^{n+1}\|}\leq e^{k_n}$$

Therefore :

\begin{eqnarray*}
 \|D(x\mapsto D_{T^{n+1}\circ \phi'_{n}(x)}T^{-(n+1)})\| & \leq & K\|D_{T^{n+1}\circ\phi'_{n}(w)}T^{-(n+1)}\|\\
 &\leq & K\|D_{T^{n+1}\circ\phi'_{n}}T^{-(n+1)}\|
 \end{eqnarray*}

\underline{\textbf{Second step :} we saturate the $n$-rectangle $\phi'_n$ in a $n+1$-rectangle $\phi'_{n+1}$ satisfying}\\
\underline{again $(H_{n+1})$.}\\

Let us denote by $\phi'_{n+1}=\phi_{n+1,y}$ the saturated $n+1$-rectangle $Sat(\phi'_n([0,1]^2))$ of the $n$-rectangle $\phi'_n$ for some $y\in Sat(\phi'_n([0,1]^2))$.  According to Corollary \ref{trois} (first item) we have $\|D(\phi_n^{'-1}\circ\phi'_{n+1})\|\leq 1$ and then by the previous step, we  get :

\begin{eqnarray*} 
\|D(x\mapsto D_{\phi'_{n+1}(x)}T^{n+1})\| & \leq & \|D(x\mapsto D_{ \phi'_{n}(x)}T^{n+1})\|\|D(\phi_n^{'-1}\circ\phi'_{n+1})\| \\
& \leq & K\|D_{\phi'_{n}}T^{n+1}\| \\
&\leq & K \|D_{\phi'_{n+1}}T^{n+1}\| 
\end{eqnarray*}

and in the same way :

\begin{eqnarray*} 
\|D(x\mapsto D_{T^{n+1}\circ \phi'_{n+1}(x)}T^{-(n+1)})\| &\leq & K \|D_{T^{n+1}\circ\phi'_{n+1}}T^{-(n+1)}\| 
\end{eqnarray*}

By the induction hypothesis (ii) and again Corollary \ref{trois} (first item) we have also

\begin{eqnarray*}
\|D\phi'_{n+1}\| & \leq & \|D\phi'_n\| \| D(\phi_n^{'-1}\circ \phi'_{n+1})\| \\
&\leq & \frac{1}{4}
\end{eqnarray*}

and similarly for all $1\leq k\leq n$ :
\begin{eqnarray*}
\|D(T^{k}\circ \phi'_{n+1})\| & \leq  & 1
\end{eqnarray*}

Obverve finally that $\phi'_n([\frac{4}{9},\frac{5}{9}]^2)\subset\phi'_{n+1}([\frac{1}{3}-K,\frac{2}{3}+K]^2)$ by Corollary \ref{trois} (second item) and that $R'_{n+1}:=\phi'_{n+1}([0,1]^2)$ satisfies $(1+K)l(\mathcal{E}_{n+1}(y)\cap R'_{n+1})\geq l(\mathcal{F}_{n+1}(y)\cap R'_{n+1}) \geq (1-K)\frac{l(\mathcal{E}_{n+1}(y)\cap R'_{n+1})}{ 
\max_{k=0,...,n}\|D_yT^k\|}$ for all $y\in R'_{n+1}$ by assumption $(G_n)$ and by Proposition \ref{encore}. \\

\underline{\textbf{Third step :} \label{pag} We cut the $n+1$-rectangle $\phi'_{n+1}$ to get a new $n+1$-rectangle $\phi_{n+1}$}\\
\underline{satisfying $\|D(T^{n+1}\circ \phi_{n+1})\|\leq 1$.}\\

\scalebox{0.45}{\includegraphics{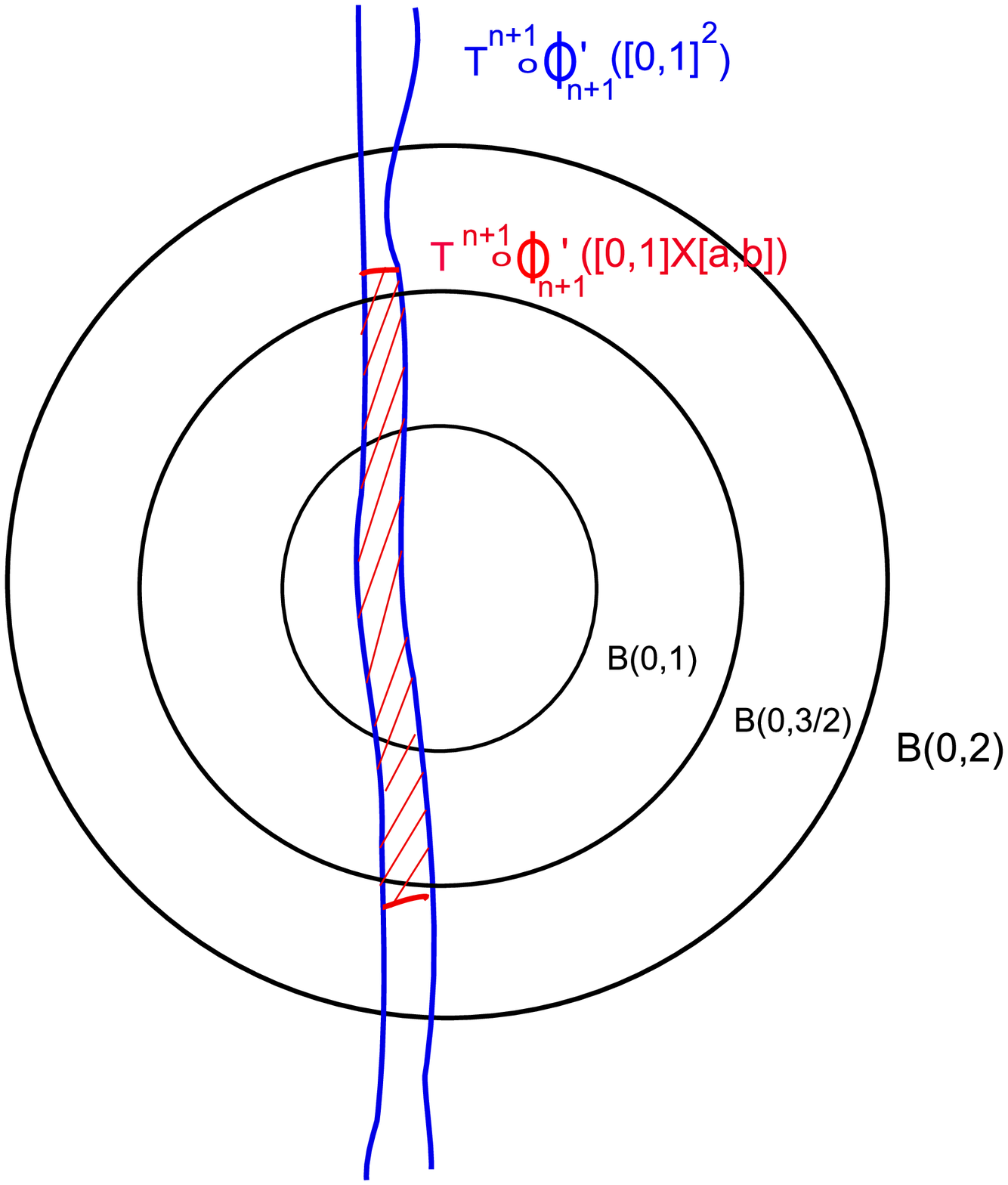}}

We consider only rectangles $\phi'_{n+1}=\phi_{n+1,y}$ whose image by $T^{n+1}$ meets the unit ball. Let $\phi'_{n+1}$ such a rectangle. We choose 
$a,b\in[0,1]$ such that the image by $T^{n+1}$ of the $n+1$-subrectangle $(\phi_{n+1},R_{n+1})$ of $(\phi'_{n+1},R'_{n+1})$, defined by $\phi_{n+1}(t,s)=\phi'_{n+1}(t,a+s(b-a))$ for all $t,s\in [0,1]$, contains $B(0,\frac{3}{2})\cap  T^{n+1}R'_{n+1}$. This is possible because $|\angle| f_{n+1}^{(n+1)}(x),f_{n+1}^{(n+1)}(y)\leq K$ according to 
Lemma \ref{angle} for all $x,y\in \phi'_{n+1}([0,1]^2)$.  Remark that $l(\mathcal{E}_{n+1}(x)\cap 
R_{n+1})=l(\mathcal{E}_{n+1}(x)\cap 
R'_{n+1})$ and $l(\mathcal{F}_{n+1}(x)\cap 
R_{n+1})\leq l(\mathcal{F}_{n+1}(x)\cap 
R'_{n+1})$ for all $x\in R_{n+1}$. Moreover either one can choose $a=0$ and $b=1$ and then $R_{n+1}=R'_{n+1}$ or $l(\mathcal{F}_{n+1}(x)\cap R_{n+1})\times \|D_xT^{n+1}\|\geq \frac{1}{2}-K$ for all $x\in R_{n+1}$. 
In both cases\footnote{Note that $l(\mathcal{E}_{n+1}(x)\cap 
R_{n+1})\leq \frac{1}{4}$ since $\|D\phi'_{n+1}\|\leq \frac{1}{4}$} we have $l(\mathcal{E}_{n+1}(x)\cap 
R_{n+1})\geq l(\mathcal{F}_{n+1}(x)\cap R_{n+1}) \geq \frac{l(\mathcal{E}_{n+1}(x)\cap R_{n+1})}{\max_{k=0,...,n+1}\|D_xT^k\|}$ for all $x\in R_{n+1}$, that is $R_{n+1}$ satisfies the assumption $(G_{n+1})$. Furthermore $\phi'_{n+1}$ satisfies assumption $(H_{n+1})$ and then so does $\phi_{n+1}$. By Corollary \ref{croco} it follows that $\|D(T^{n+1}\circ \phi_{n+1})\|\leq 2+K$.
Finally by distinguishing the cases $0<,=a\leq b<,=1$ it is easily seen that 
$\phi'_{n+1} ([\frac{1}{3}-K,\frac{2}{3}+K]^2)\cap B(n+2,1)\subset \phi_{n+1}([\frac{1}{3}-K,\frac{2}{3}+K]\times [\frac{1}{6}-K,\frac{5}{6}+K])$.\\ 

In summary we have :
\begin{itemize}
\item  $\mathcal{H}^n_{\mathcal{T}}(\chi^+,\chi^-,\gamma,C)\cap \mathcal{H}(\mathcal{K}_{n-1})\cap B(n+1,1)\subset \bigcup_{\phi_n\in \mathcal{F}_n}\phi([\frac{4}{9},\frac{5}{9}]^2)$ by induction hypothesis (vi) ;
\item $\mathcal{H}^{n+1}_{\mathcal{T}}(\chi^+,\chi^-,\gamma,C)\cap\mathcal{H}(\mathcal{K}_n)\cap \left(\bigcup_{\phi_n}\phi_n([\frac{4}{9},\frac{5}{9}]^2)\right)\subset \bigcup_{\phi'_n}\phi'_{n}([\frac{4}{9},\frac{5}{9}]^2)$ according to the first step ;
\item $\phi'_n([\frac{4}{9},\frac{5}{9}]^2)\subset\phi'_{n+1}([\frac{1}{3}-K,\frac{2}{3}+K]^2)$ according to the  second step ;
\item $\phi'_{n+1}([\frac{1}{3}-K,\frac{2}{3}+K])\cap B(n+2,1)\subset \phi_{n+1}([\frac{1}{3}-K,\frac{2}{3}+K]\times [\frac{1}{6}-K,\frac{5}{6}+K])$ according to the third step. 
\end{itemize}

Therefore by dividing the $n+1$-rectangles $(\phi_{n+1},R_{n+1})$ into at most $20$ subrectangles and by reparametrizing them by affine contractions we can ensure that 
$\mathcal{H}^{n+1}_{\mathcal{T}}(\chi_+,\chi_-,\gamma,C)\cap \mathcal{H}(\mathcal{K}_n)\cap B(n+2,1)\subset \bigcup_{\phi_{n+1}}\phi_{n+1}([\frac{4}{9},\frac{5}{9}]^2)$.

\end{demo}

\section{Beyond $\mathcal{C}^2$ surface diffeomorphisms}

\subsection{Non-invertible maps}
In the previous proof of Proposition \ref{main} one only needs to estimate the entropy locally and therefore only local invertibility is required. Theorem \ref{gg} and therefore Theorem \ref{super}  are also valid in the context of  surface local diffeomorphisms. Indeed in the proof of Theorem \ref{gg} global invertiblity was just used to prove by Ruelle-Margulis inequality  that ergodic measures with positive Newhouse local entropy and therefore with positive entropy have one positive and one negative Lyapunov exponent. This fact is still true for surface local diffeomorphisms. Indeed the author proved in \cite{superbur} the following "local invertible" version of Ruelle-Margulis inequality :

\begin{lem}\label{ruel}\cite{superbur}
Let $M$ be a compact surface and let $T:M\rightarrow M$ be a local diffeomorphism. There exists $\epsilon>0$ such that for any ergodic measure $\nu$, we have :
$$h^{New}(M|\nu,\epsilon) \leq \min(\chi^+_0(\nu),-\chi_0^-(\nu))$$
\end{lem}

 The lemma remains true in higher dimensions by replacing $\chi^+_0$ by the sum of the positive Lyapunov exponents and $\chi^-_0$ by the sum of the negative Lyapunov exponents \cite{superbur}.\\
 
 The conjecture \ref{mainconj} is still open for $\mathcal{C}^2$ noninvertible surface maps with critical points.

\subsection{Higher regularity}
In Yomdin's theory \cite{Yoma},\cite{Yomb},\cite{Gr},\cite{BurY} one reparametrizes the Bowen ball by semi-algebraic contracting maps. Here we avoid 
the semi-algebraic tools involving in this 
theory but we are then reduced to the $\mathcal{C}^2$ case. Assume $T:M\rightarrow M$ is a $\mathcal{C}^r$ map with $r\in \N$ and $r>2$ and $\phi_n:[0,1]^2\rightarrow M$ is a $n$-rectangle. To bound the oscillation of the
 derivative $D_{\phi_n(x)}T^n$ one can try to control the $r$ derivative $D^r(x\mapsto D_{\phi_n(x)}T^n)$ and then  approximate the derivative at one point  by its Lagrangian polynomial as in Yomdin's theory. Following the proof of  Proposition \ref{main} the number of subdivisions at step $n$ we need to bound
$\|D^r(x\mapsto D_{\phi_n(x)}T^n)\|$ by $\|D_{\phi_n}T^n\|$ is of order $Ke^{\frac{2k_n}{r-1}}$ which is consistent with Conjecture \ref{mainconj}. 
However to bound the dynamical complexity of the Lagrangian polynomial of $D_{\phi_n(x)}T^n$ as in Yomdin's theory we have to compose $\phi_n$ with a semi-algebraic map $\psi:[0,1]^2\rightarrow [0,1]^2$ and no more a homothety. Then it seems difficult to relate the new map $ \phi'_n=\phi_n\circ \psi$ with a $n$- or  $n+1$-rectangle.

\subsection{Higher dimensions}
We used in the proof typical one dimensional arguments to analyse the finite time stable and unstable manifolds (in  particular the integrability of finite time stable and unstable fields). Moreover ergodic measures with positive entropy need not to be hyperbolic (i.e. with nonzero Lyapunov exponents) in higher dimensions. Finally  the creation of horseshoes is not the only mechanism to create entropy. Indeed a horseshoe generates many periodic points but M.Herman \cite{her} built a minimal diffeomorphism on a compact 4-dimensional manifold with positive topological entropy.

\it{E-mail  address :}  \texttt{burguet@math.polytechnique.fr}

\end{document}